\documentclass[12pt]{article}
\usepackage[final]{epsfig}
\usepackage{graphics}
\usepackage{amsmath}
\usepackage{amsfonts}
\usepackage{latexsym}
\usepackage{amssymb}
\usepackage{graphicx}
\usepackage{url}
\usepackage{nopageno}
\usepackage[all]{xy}
\usepackage{epstopdf}

\newtheorem{lemma}{Lemma}[section]
\newtheorem{proposition}[lemma]{Proposition}

\newtheorem{theorem}{Theorem}
\newtheorem{definition}[lemma]{Definition}
\newtheorem{corollary}[lemma]{Corollary}

\newcommand{\proofend}{$\Box$\bigskip}

\def\proof{\paragraph{Proof.}}

\begin{document}

\title{Discrete Conics}

\author{Emmanuel Tsukerman\footnote{Stanford University, emantsuk@stanford.edu}
}

\date{}
\maketitle

\begin{abstract}

In this paper, we introduce {\it discrete conics}, polygonal analogues of conics. We show that discrete conics satisfy a number of nice properties analogous to those of conics, and arise naturally from several constructions, including the discrete negative pedal construction and an action of a group acting on a focus-sharing pencil of conics.

\end{abstract}

\section{Introduction} \label{intro}

In \cite{Discrete Parabolas}, we showed that a certain family of
Simson polygons, polygons which admit a point whose projections
into the sides of the polygon are collinear, can be fruitfully viewed
as discrete analogs of the parabola. In this paper, we find a family
of polygons which are discrete analogs of a general conic. 

\begin{figure}[hbtp] \label{fig: discConic}
\centering
\includegraphics[width=3.4in]{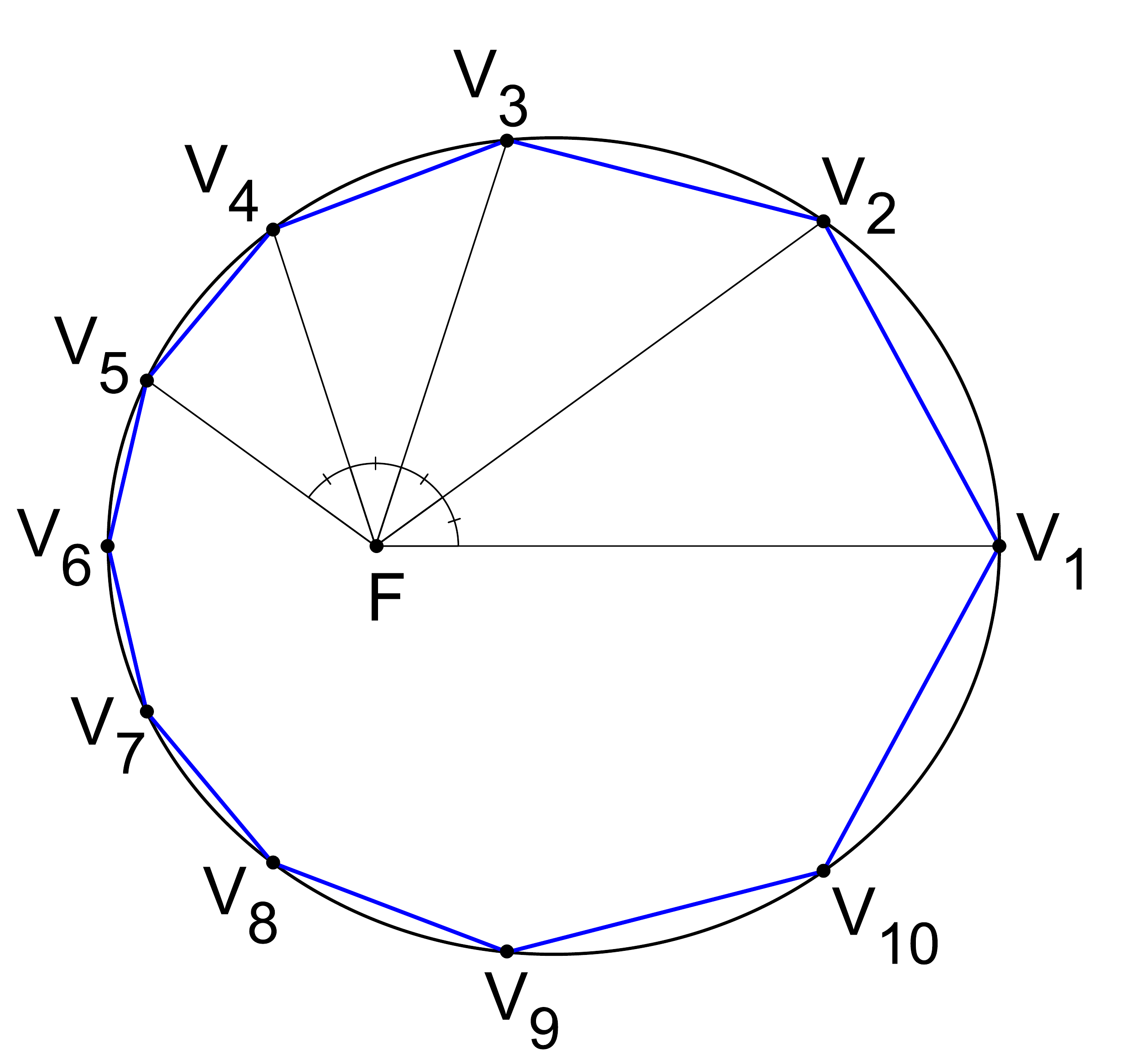}
\caption{A discrete conic.}
\end{figure}

Intuitively, a regular polygon is the best candidate for being called a "discrete circle" - it has the maximal amount of symmetry as well as other geometric properties that are similar to those of a circle. In the same spirit, we define discrete analogues of conics as follows (see Figure \ref{fig: discConic}):

\begin{definition} (Discrete Conic) Let $C$ be a conic with focus $F$. A discrete conic
is a polygon $V_{1}V_{2} \cdots V_{n}$ such
that for every $i\in {1,2,\ldots,n}$ and some fixed $\theta \in \mathbb{R}$, $\angle V_{i}FV_{i+1}=\theta$.\footnote{When $C$ is a hyperbolas, parabola or a degenerate conic, we allow the rays used in measuring an angle to wrap around the plane.}
\end{definition}

We will show later that these discrete conics are projective images of regular polygons. One of our goals in this paper will be to show how these discrete conics arise naturally in two constructions: via a group acting on a pencil of focus-sharing conics and through the discretized negative pedal construction.
In the last section, we will prove a number of nice properties of discrete conics.

\section{Properties of Discrete Conics}

We now list some of the nice properties of discrete conics. Proofs can be found in section \ref{Proofs}.

\begin{theorem}

\label{Thm: tangency points projectively regular}A discrete conic
is projectively equivalent to a regular polygon.
\end{theorem}

Recall that a {\it Poncelet Polygon} is a polygon which circumscribes a conic and is inscribed in a conic.

\begin{theorem} \label{Thm: Poncelet Polygon}A discrete conic
is a Poncelet Polygon and the circumscribing and inscribed conics share a focus.
\end{theorem}

Let us fix notation. The discrete
conic $D=V_{1}V_{2}\cdots V_{n}$ lies on conic $C$ with focus
$F$. The vertices of $D$ satisfy $\angle V_{i}FV_{i+1}=\theta$.
We set $S_{i}=V_{i}V_{i+1}$. The sides $S_{i}$ are tangent to a conic with focus $F$ which we
will denote by $C_{L}$. Let $F'$ be its other focus. We define $F$
and $F'$ to be the \textit{focii of the discrete conic} $D$. Finally,
we set $M_{i}=V_{i}V_{i+1}\cap C_{L}$ to be a tangency point of the
inner conic.

\begin{theorem}

\label{thm:Diagonals coincide at P'}Let $D=V_{1}\cdot\cdot\cdot V_{n}$
be a discrete conic with $n$ sides. If $n$ is even then the diagonals
$V_{j}V_{j+\frac{n}{2}}$ of $D$ are coincident at $F$. If $n$
is odd then chords $V_{j}M_{j+\frac{n-1}{2}}$ are coincident at $F$.
\end{theorem}

An analog of the reflective property of a conic
exists in the case when $D$ is closed. Unlike the previous results, this is not a projective property and does not follow from Theorem \ref{Thm: tangency points projectively regular}. We will symbolically represent
a path taken by a ray via a sequence of alternating positions and
arrows. For instance, if $r$ is a ray which starts at $F'$, reflects
off a side $S$ and then passes through $F$, then we will say that
$r$ traverses the path $F'\rightarrow S\rightarrow F$. 

\begin{figure}[hbtp]
\centering
\includegraphics[width=2.4in]{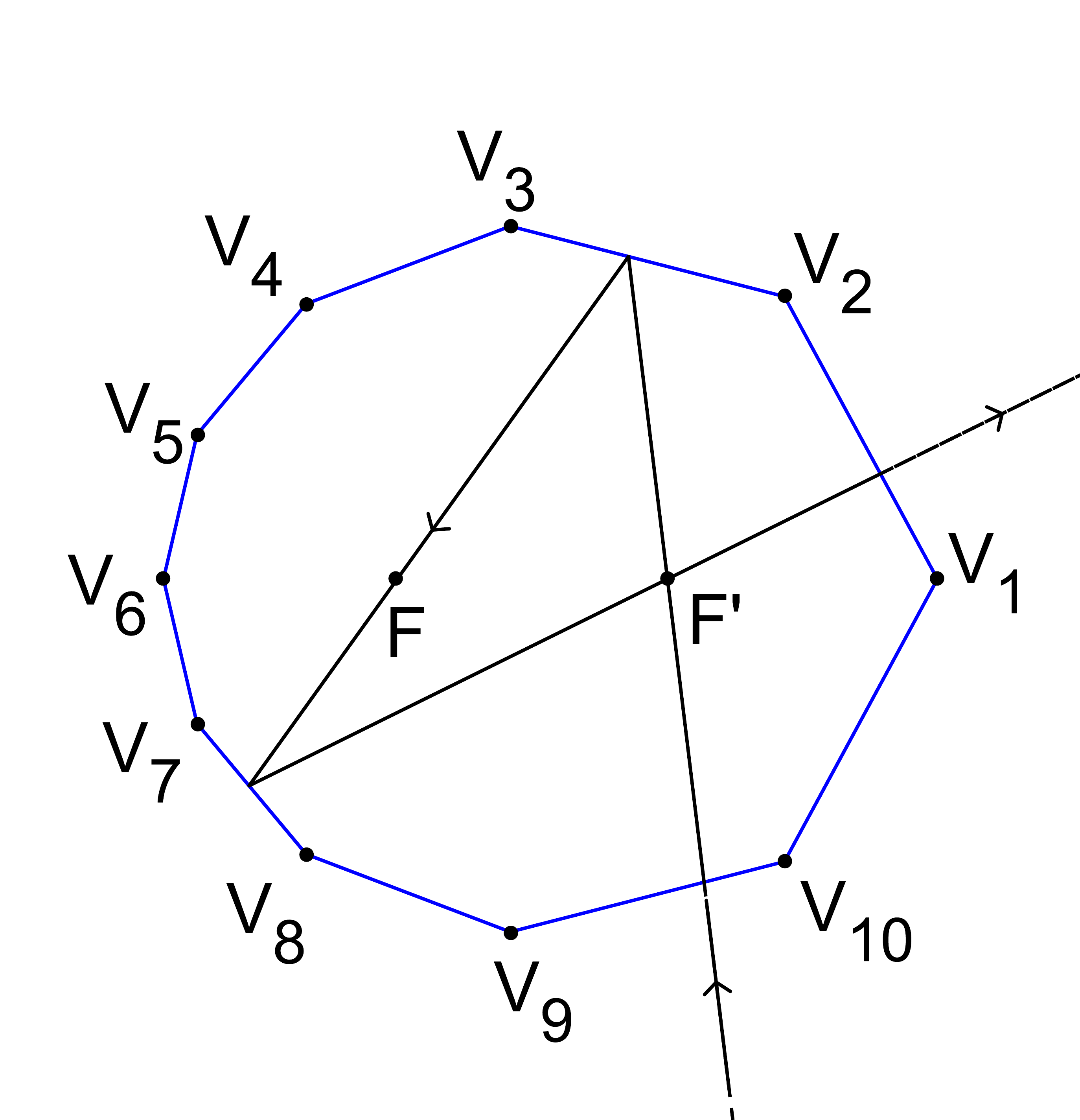}
\includegraphics[width=2.4in]{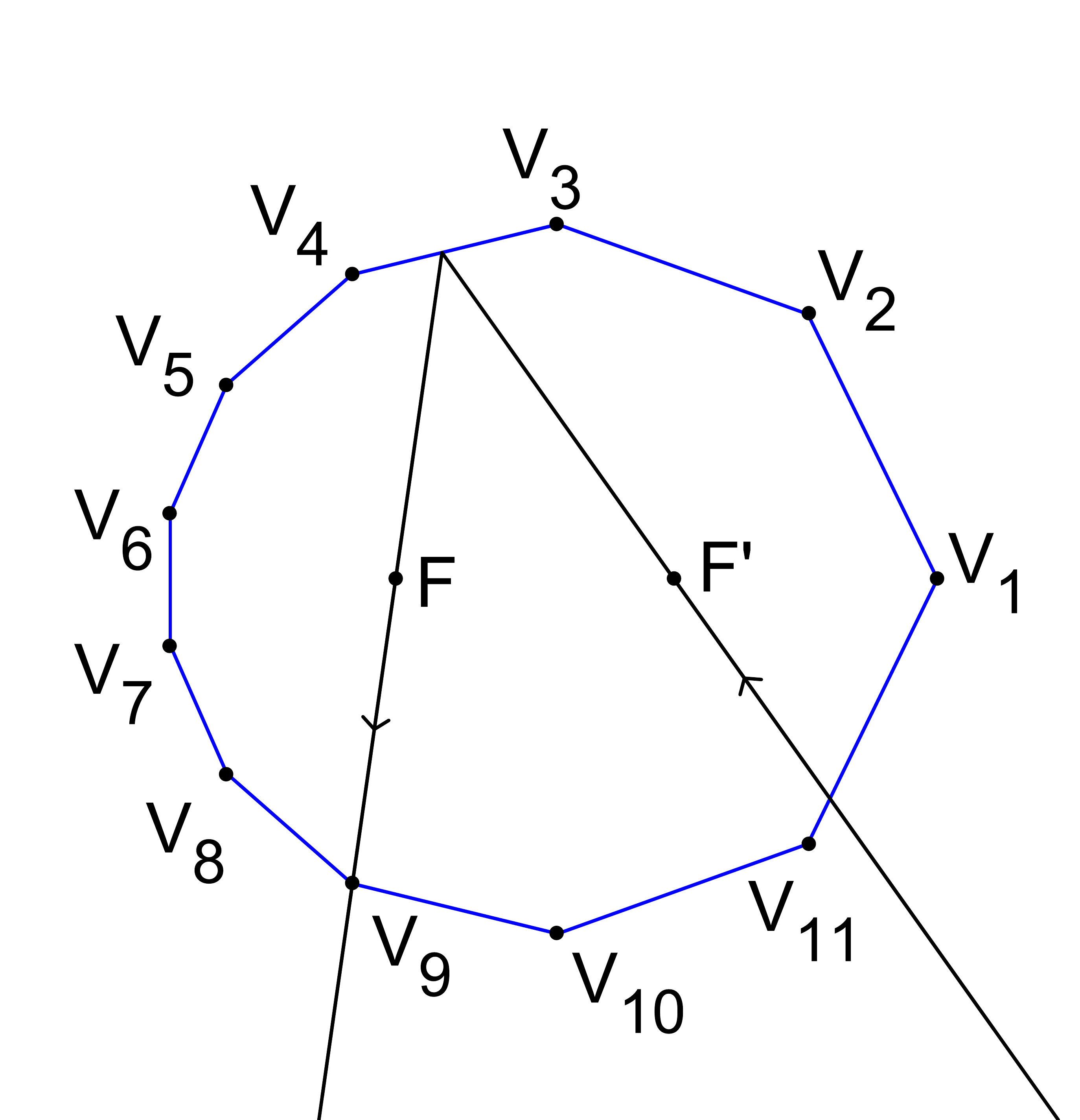}
\caption{Optic (reflective) property for even- and odd-sided discrete conics. }
\end{figure}

\begin{theorem}

\label{thm:Reflective Property}(Reflective Property) Let $r$ be
a ray which traverses the path $F'\rightarrow S_{j}\rightarrow F$.
If the discrete conic $D$ is an even-sided polygon then $r$ traverses
\textup{
\[
F'\rightarrow S_{j}\rightarrow F\rightarrow S_{j+\frac{n}{2}}\rightarrow F'.
\]
}

\textup{If $n$ is odd-sided then $r$ traverses the path 
\[
F'\rightarrow S_{j}\rightarrow F\rightarrow V_{j+\frac{n-1}{2}}.
\]
}
\end{theorem}

As a particular example of Theorem \ref{thm:Reflective Property},
consider the case when $C$ is a circle. Then $D$ is a regular polygon.
In this case, the Theorem states that if the polygon is even-sided
and regular then a ray through the center which reflects back to the
center must consequently reflects back to the center again. In the
case that the regular polygon is odd-sided then a ray passing through
the center that reflects back to the center must subsequently pass
through one of the vertices.

We also have an analog of the Isogonal Property of a conic. Recall
that if tangents at $X$ and $Y$ to a conic intersect at $Z$ and
$F,F'$ are the focii of the conic, then $\angle FZX=\angle F'ZY$
\cite{Akopyan Geometry of Conics}. Analogously,

\begin{theorem} \label{thm:isogonal property}

(Isogonal Property) If $Z$ is a point lying on the intersection of
sides $S_{i}$ and $S_{j}$ of two sides of the discrete conic $D$
then $FZ$ bisects angles $\angle M_{i}FM_{j}$, $\angle V_{i}FV_{j+1}$
and $\angle V_{i+1}FV_{j}$. Moreover, $FZ$ and $F'Z$ are isogonal
with respect to $S_{i}$ and $S_{j}$.
\end{theorem}

\begin{figure}[hbtp]
\centering
\includegraphics[width=2.5in]{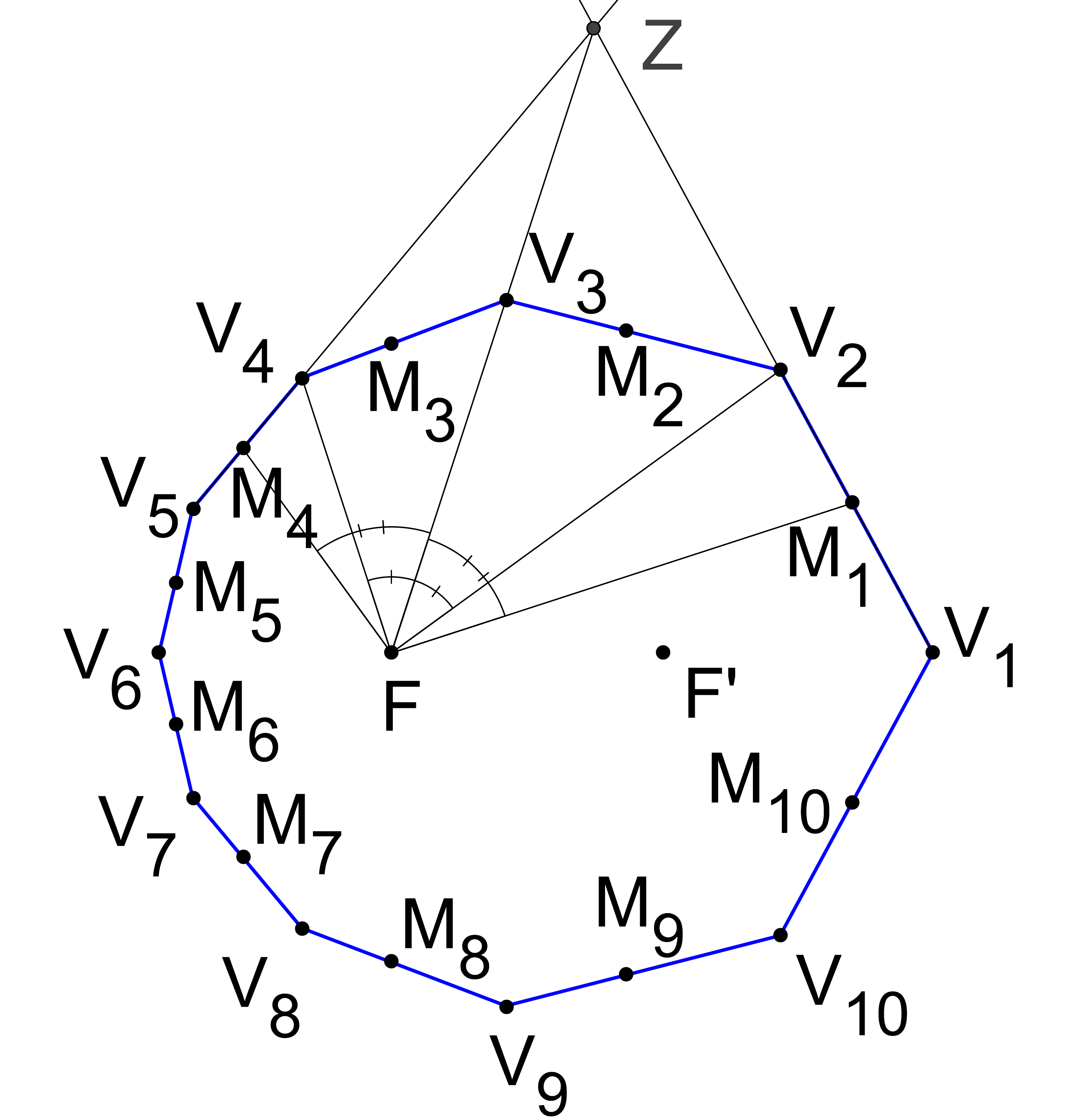}
\includegraphics[width=2.5in]{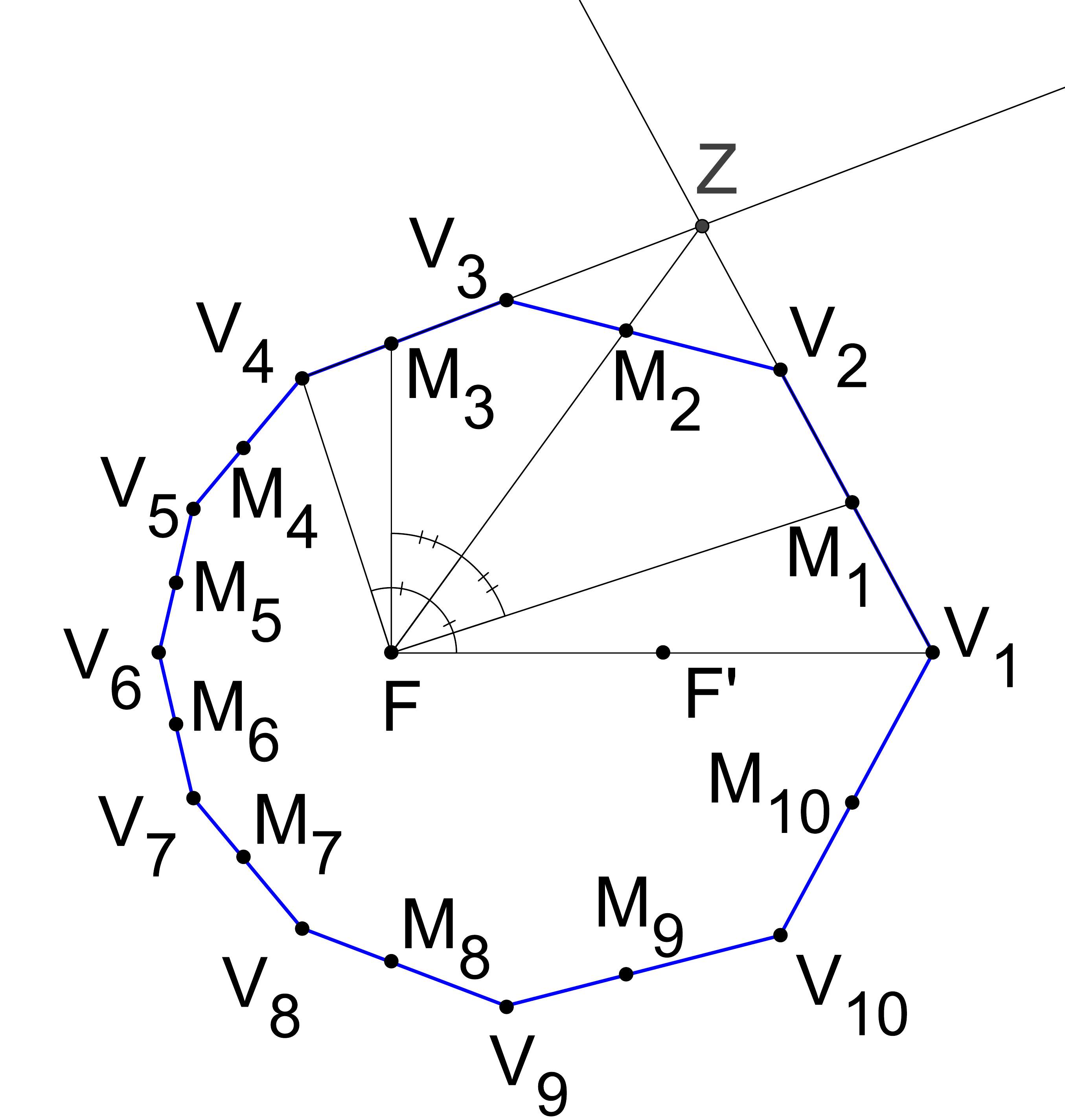}
\caption{Isogonal properties of discrete conics. }
\end{figure}

\begin{figure}[hbtp]
\centering
\includegraphics[width=5.0in]{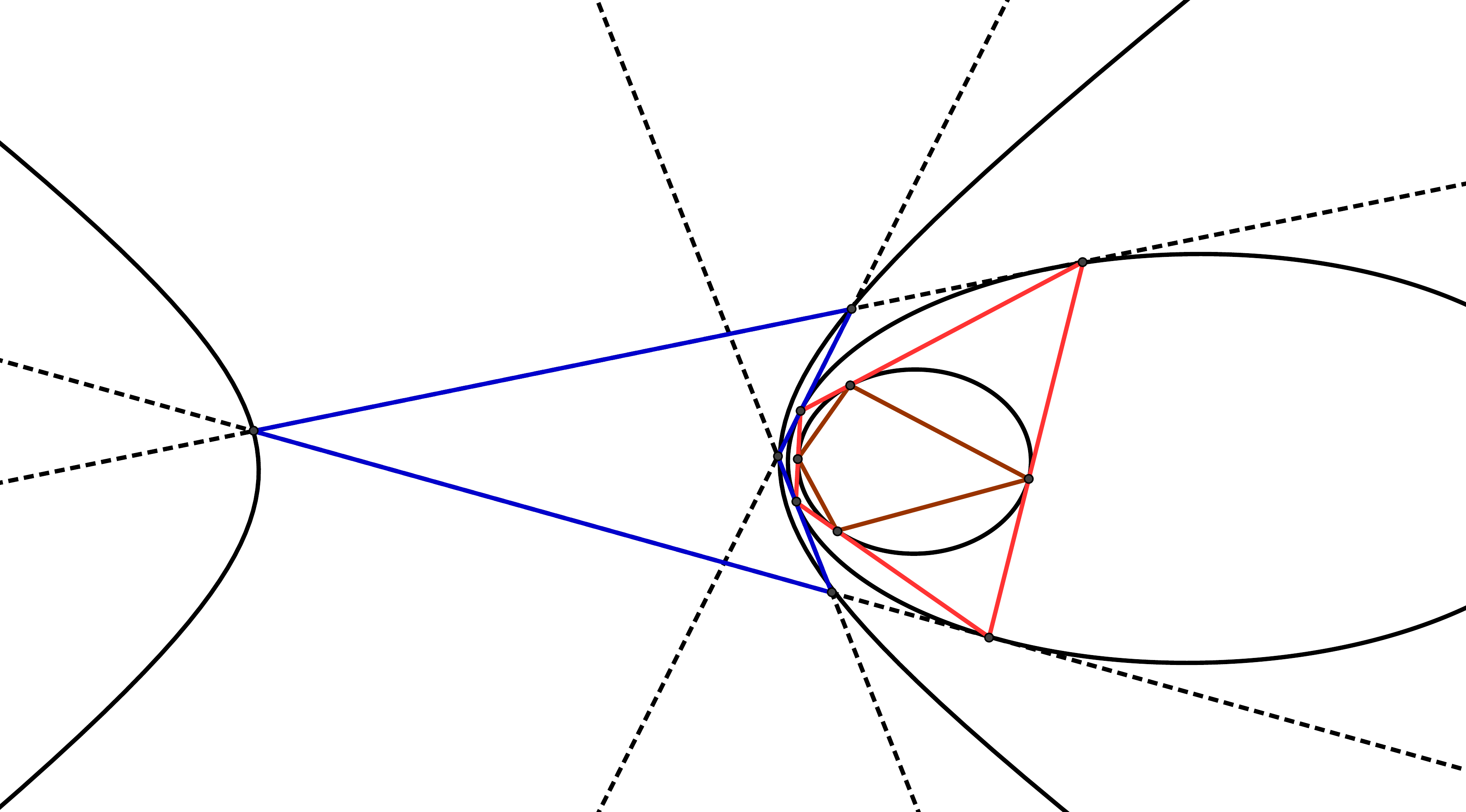}
\caption{\label{fig:Pencil Discrete Conic}Three members of an infinite family of discrete conics. Given a discrete
conic $D_{i}$ on conic $C_{i}$, the subsequent discrete conic $D_{i+1}$
has vertices which are the intersections of tangents of $C_{i}$ at
the vertices of $D_{i}$. All conics share a single focus.}
\end{figure}

We can describe a discrete conic by the conic it is on (parameter
$t$), the angle $\theta$ of the relation $\theta=\angle V_{i}FV_{i+1}$
and the phase angle $\phi\mod\theta$. Thus we will write $D_{t,\theta,\phi}$.
The group operations that can be applied to this discrete conic correspond to multiples of angle $\theta$. In particular, $G_{k\theta}$ sends $D_{t,\theta,\phi}$ to
\[
G_{k\theta}(D_{t,\theta,\phi})=D_{t\sec^{2}(\frac{k\theta}{2}),\theta,\phi+\frac{k\theta}{2}}.
\]

Using this fact, we will discuss the relation between discrete conics and the Poncelet Grid Theorem.
R. Schwartz's Poncelet Grid Theorem states that if $P=P_1 P_2\cdots P_{n}$ is a Poncelet polygon then the intersection points of lines $L_i=P_i P_{i+1}$ and $L_{i+k}=P_{i+k} P_{i+k+1}$ lie on a conic $c'$ and form a polygon $P'$ projectively equivalent to $P$ (see \cite{Poncelet Grid} and \cite{Tabachnikov and Levi}). More can be said when $P$ is a discrete conic.

\begin{theorem} \label{thm:Poncelet Grid} If $P$ is a discrete conic lying on conic $c$ then $P'=G_{k\theta}(H_{\theta}(P))$. In particular, $P'$ is a discrete conic having the same angle parameter $\theta$ and lying on a a conic which shares a focus with $c$. 
\end{theorem}

A consequence of what we have shown
is that there exists an infinite family of discrete conics, each discrete conic
lying on a conic of this pencil, and with tangents to the corresponding conic at the vertices of the discrete conic
forming the sides of the subsequent discrete conic (see Figure \ref{fig:Pencil Discrete Conic}).

\begin{figure}[hbtp]
\centering
\includegraphics[width=4in]{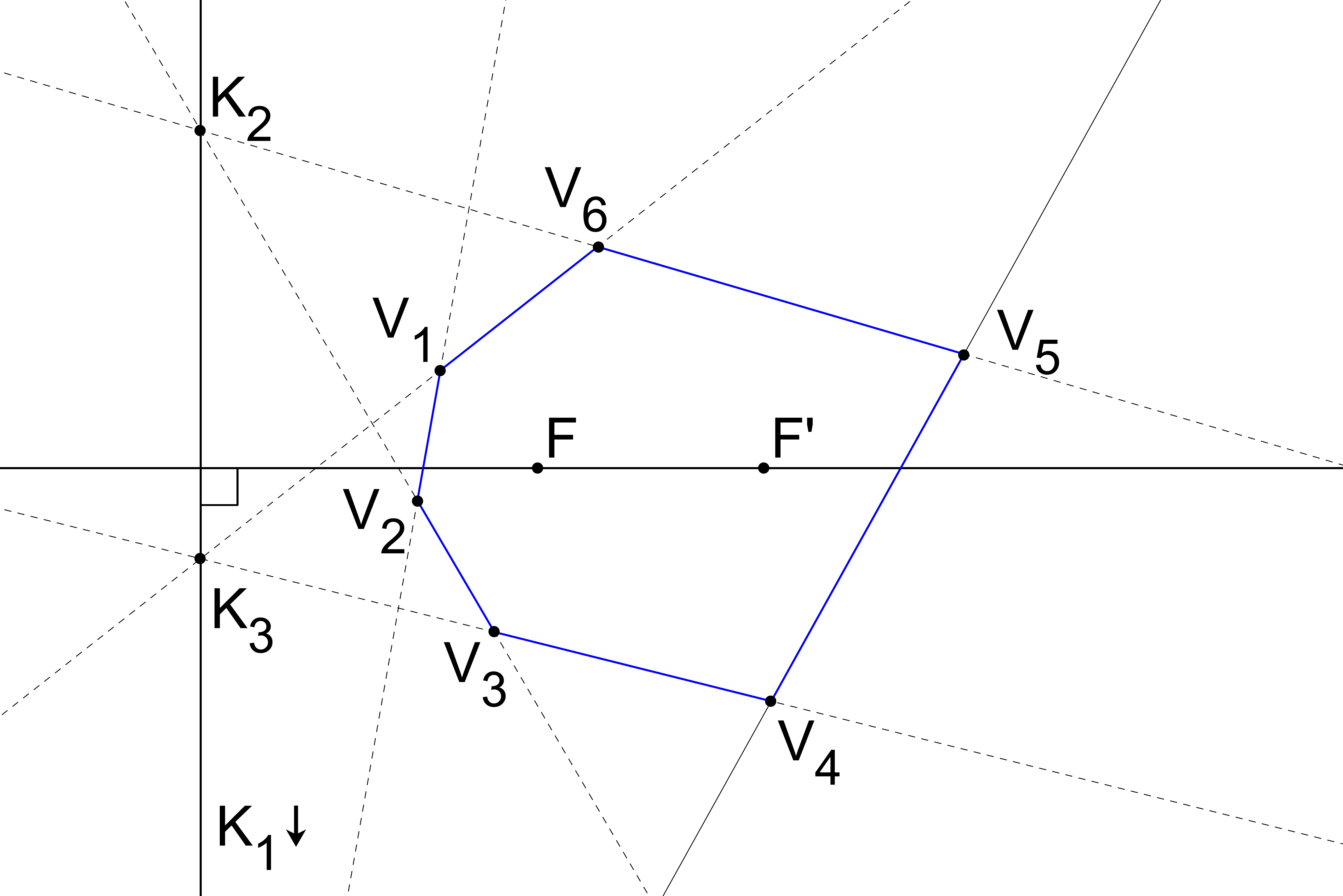}
\caption{Theorem \ref{thm: pascal} for a hexagon. We set $K_i=V_i V_{i+1} \cap V_{i+n}V_{i+n+1}$ for $i=1,2,...,2n-1$ to be the intersections of opposite sides. These are collinear and their line is orthogonal to $FF'$. Moreover, these form a discrete conic, so that the angles as seen from $F$ are equal.}
\end{figure}

We conclude with an interesting result, reminiscent of Pascal's famous hexagon theorem. 

\begin{theorem} \label{thm: pascal} Let $D=V_1V_2 \cdots V_{2n}$ be an even-sided discrete conic with focii $F,F'$. The intersections of opposite sides lie on a line orthogonal to $FF'$ and form another discrete conic. In the special case when $C$ is a parabola, this line is the directrix of the parabola.
\end{theorem}

\section{A Group Acting on a Pencil of Focus-Sharing Conics} \label{sec:Group-Acting-on}

Let $C$ be a circle with center $O$. Let $X$ and $Y$ be variable
points on $C$. Define $G_{\theta}(C)$ to be the locus of points
$Z$ such that $Z$ is the intersection of tangent to $C$ at $X$
and $Y$ for $X$ and $Y$ satisfying $\angle XOY=\theta$. We define
$H_{\theta}(C)$ to be the locus enveloped by lines $z=XY$ for $X$
and $Y$ satisfying $\angle XOY=\theta$. We restrict $\theta\in[0,\pi)$.

\begin{figure}[hbtp]
\centering
\includegraphics[width=2.6in]{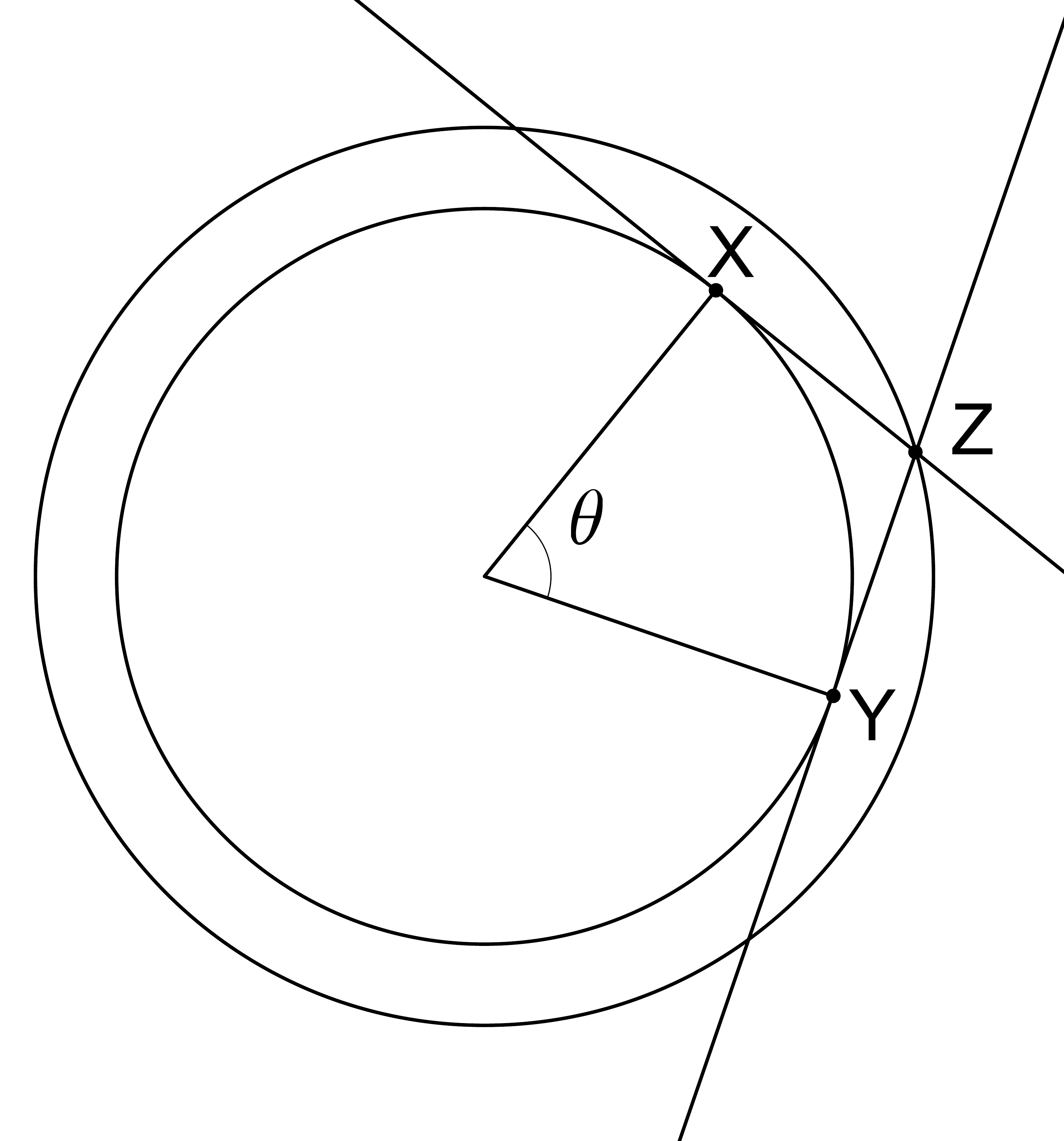}
\includegraphics[width=2.6in]{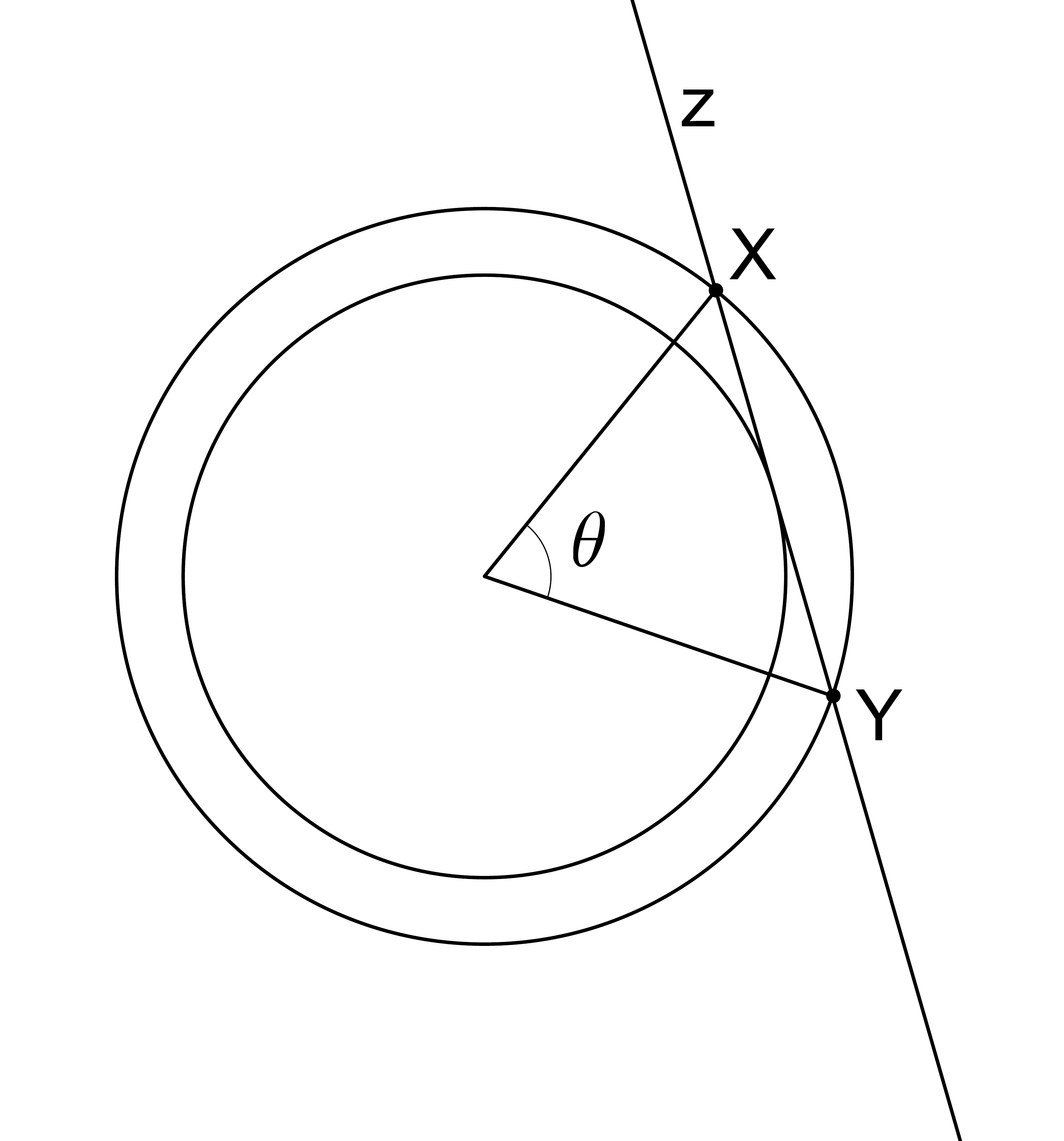}
\caption{The operations $G_{\theta}(C)$ and $H_{\theta}(C)$.}
\end{figure}

\begin{lemma} \label{lem:Group on concentric circles} Consider a family of concentric plane-foliating circles $OF$. The operations $G$ and $H$ define an abelian group acting on $OF$ with $G_{\theta}=H_{\theta}^{-1}$ and $G_{\theta}\circ G_{\phi}=G_{2\cos^{-1}(\cos(\frac{\theta}{2})\cos(\frac{\phi}{2}))}$.
\end{lemma}

\proof
Let $R$ be the radius of circle $C$ with center $O$. From right
triangle $\triangle OXZ$, we see that the radius of $G_{\theta}(C)$
is $R\sec(\frac{\theta}{2})$. On the other hand, let $M$ be the
midpoint of $XY$. From right triangle $\triangle OMX$, the radius
$OM$ of $H_{\theta}(C)$ is equal to $R\cos(\frac{\theta}{2})$.
It follows that $G_{\theta}=H_{\theta}^{-1}$ and $G_{\theta}\circ G_{\phi}=G_{2\cos^{-1}(\cos(\frac{\theta}{2})\cos(\frac{\phi}{2}))}$.
The identity is just $G_{0}$.
\proofend

The map $\theta\rightarrow\cos\frac{\theta}{2}$ sending angle $\theta\in[0,\pi)$
to the number $\cos\frac{\theta}{2}$ shows that the group is isomorphic
to the positive real numbers under multiplication.

We now extend this group operation to conics. Consider a conic with
focus $F$. Define $G'_{\theta}(C)$ to be the locus of points $Z$
such that $Z$ is the intersection of tangents at $X$ and $Y$ satisfying
$\angle XFY=\theta$. We define $H'_{\theta}(C)$ to be the locus
enveloped by lines $z$ such $z=XY$ for $X$ and $Y$ satisfying
$\angle XFY=\theta$.

\begin{figure}[hbtp]
\centering
\includegraphics[width=2.6in]{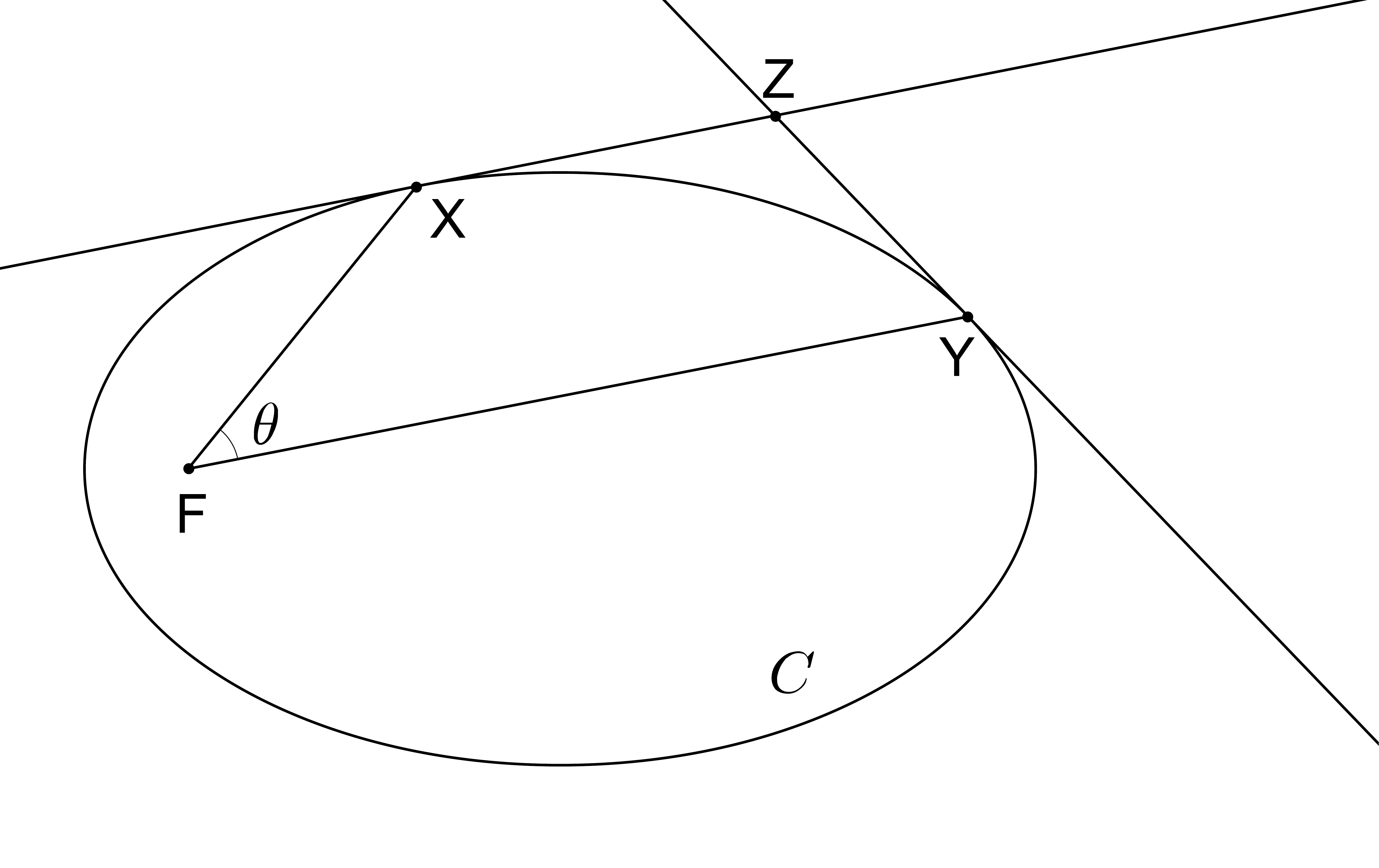}
\includegraphics[width=2.6in]{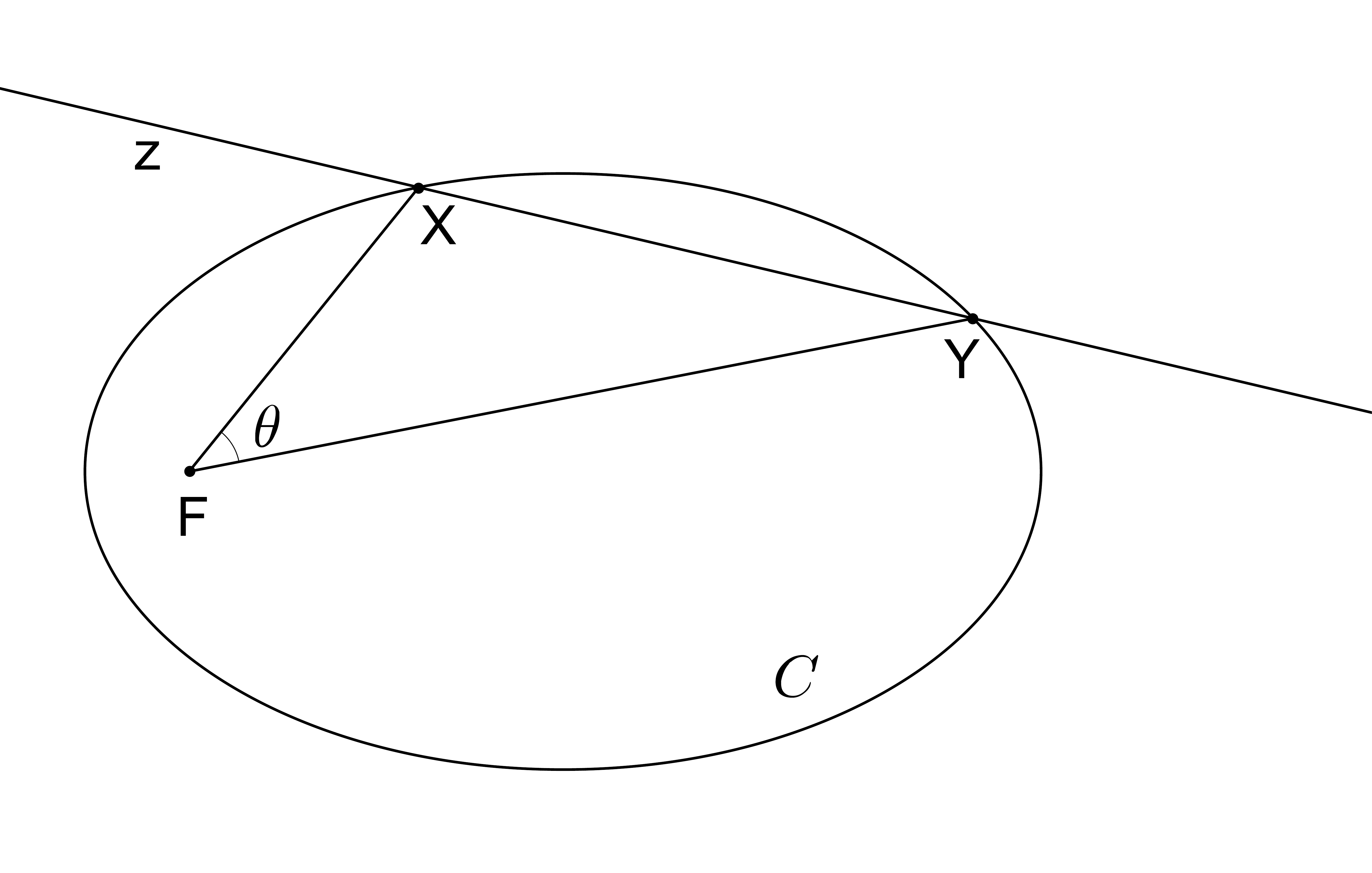}
\caption{The operations $G'_{\theta}(C)$ and $H'_{\theta}(C)$.}
\end{figure}

\begin{theorem} \label{thm:conic locus} 
Let $C$ be a conic with focus $F$. Then
$G'_{\theta}(C)$ and $H'_{\theta}(C)$ are conics with focus $F$.
\end{theorem}

For the proof, we will use a technique suggested by Arseniy Akopyan, namely reciprocation (also known as polar duality). For a detailed discussion and proofs we refer the reader to \cite{Geometry Revisited}. Consider a circle with center $O$.  Each point $P$ other than $O$ determines a corresponding line $p$, the {\it polar} of P, which is the line perpendicular to $OP$ and passing through the inverse of $P$.  Conversely, given a line $p$, we call the point which is the inverse of the foot of the perpendicular from $O$ to $p$  its {\it pole}. The transformation taking a set of points to their polar lines or a set of lines to their poles is called reciprocation. We do not distinguish the points of a curve from the envelope of tangents to the curve.  \\
A fact we will employ in the proof is that a conic is the reciprocal of a circle and vice verse. Moreover, point $O$ is the focus of the conic. To see why this is true, consider the reciprocal of the tangent lines to the conic. To reciprocate a line, we first find the foot of the perpendicular from $O$. Since the pedal of a conic with respect to its focus is a circle, the locus of feet of perpendiculars is the pedal circle. We then invert these to obtain the points of the reciprocal of the conic. Hence these form a circle, and in particular are the inverse image of the pedal circle.
  \\

\begin{figure}[hbtp]
\centering
\hspace{-50pt}
\includegraphics[width=2.8in]{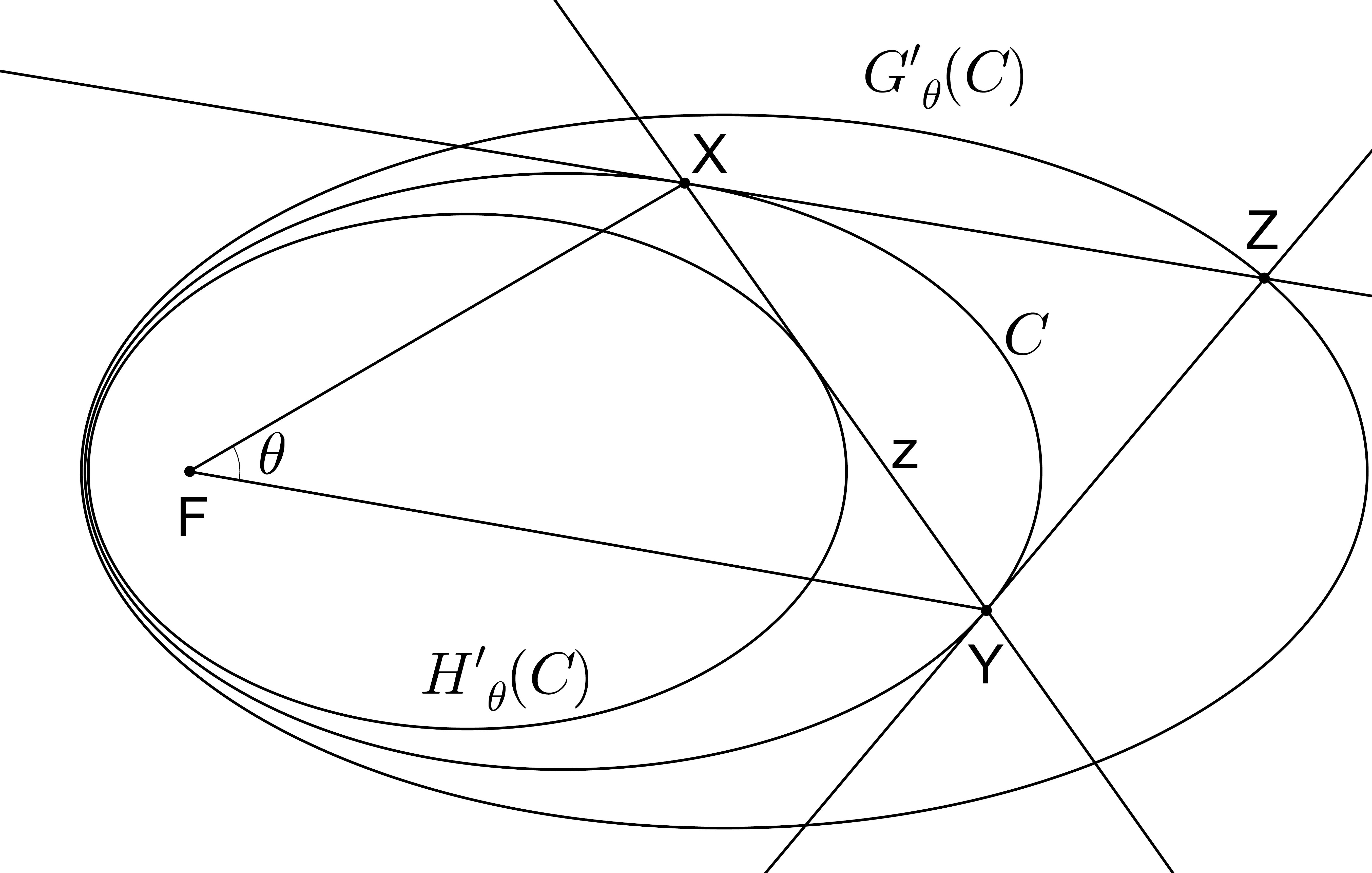}
\caption{The operations $G'_{\theta}(C)$ and $H'_{\theta}(C)$ produce conics
which share focus $F$.}
\end{figure}

\noindent{\bf Proof of Theorem \ref{thm:conic locus}.}
Let $x$ and $y$ be the tangent lines at $X$ and $Y$, respectively.
We apply a duality transformation $d$ about a circle centered at
$F$. The transformation takes the conic $C$ to a circle $d(C)$
upon which points $d(x)$ and $d(y)$ lie. It is not hard to see that we may choose the circle of inversion so that $F$ is inside of $d(C)$. Point $Z$ maps to the
line $d(Z)$ passing through $d(x)$ and $d(y)$. The images $d(X)$
and $d(Y)$ are the tangents to $d(C)$ at $d(x)$ and $d(y)$. Let $X'$ and $Y'$ be the feet of the perpendiculars $Fd(X)$ and $Fd(Y)$. Since $\angle XFY=\theta$, we have $\angle X'FY'=\theta$. Let $\bar c$ be the center of $d(C)$. Then 
\[
\angle d(x)\bar c d(y)=\pi-\angle d(x)d(z)d(y)=\pi-\angle X' d(z) Y' = \angle X'FY'=\theta
\]
Therefore the envelope of $d(Z)$ is a
concentric circle.  Applying $d$ to this concentric circle shows
that $G'_{\theta}(C)$ is a conic with focus $F$. Similar reasoning
shows that $H'_{\theta}(C)$ is also a conic with focus $F$. 
\proofend

\begin{corollary} \label{cor:Commutative Diagrams}
Let $C$ be a conic with focus $F$,
and $d$ be the duality transform about $F$. The following diagrams
commute for every $\theta$:

\begin{equation*}  
\xymatrix@R+2em@C+2em{
  C \ar[r]^-{d} \ar[d]_-{H'_\theta} & d(C) \ar[d]^-{G_\theta} \\
  H'_\theta(C) \ar[r]_-{d}& d(H'_\theta(C))=G_\theta(d(C))   } 
\end{equation*}

\begin{equation*}  
\xymatrix@R+2em@C+2em{
  C \ar[r]^-{d} \ar[d]_-{G'_\theta} & d(C) \ar[d]^-{H_\theta} \\
  G'_\theta(C) \ar[r]_-{d}& d(G'_\theta(C))=H_\theta(d(C))   } 
\end{equation*}
\end{corollary}

\proof
The commutativity of the diagrams follows from the proof of Theorem
\ref{thm:conic locus}.
\proofend

Using the result of Corollary \ref{cor:Commutative Diagrams}, we
can show that $G'$ and $H'$ form a group isomorphic to that of $G$
and $H$.
\begin{theorem}

The operations $G'$ and $H'$ define an abelian group with $G_{\theta}' =  {H'_{\theta}}^{-1}$ and $G_{\theta}'\circ  G_{\phi}'=G_{2\cos^{-1}(\cos(\frac{\theta}{2})\cos(\frac{\phi}{2}))}'$.
\end{theorem}
\proof
By Corollary \ref{cor:Commutative Diagrams}, we have $H_{\theta}'=d\circ G_{\theta}\circ d^{-1}$
and $G_{\theta}'=d\circ H_{\theta}\circ d^{-1}$. Therefore $G_{\theta}'\circ H_{\theta}'=d\circ H_{\theta}\circ d^{-1}\circ d\circ G_{\theta}\circ d^{-1}=I$
and similarly $H_{\theta}'\circ G_{\theta}'=I$. We also have $G_{\theta}'\circ G_{\phi}'=d\circ G_{\theta}\circ d^{-1}\circ d\circ G_{\phi}\circ d^{-1}=d\circ G_{\theta}\circ G_{\phi}\circ d^{-1}$.
Finally, $G_{\theta}'\circ G_{0}'=d\circ G_{\theta}\circ G_{0}\circ d^{-1}=d\circ G_{\theta}\circ d^{-1}=G_{\theta}'$.

Consider the map $\Theta:G(CF)\rightarrow G(OF)$ from the group acting
on conics to the group acting on the circles given by 
\[
\Theta(G_{\theta}')=d^{-1}\circ G_{\theta}'\circ d=H_{\theta}.
\]

This map is clearly a bijection. Moreover, 
\[
\Theta(G_{\theta}'\circ G_{\phi}')=d^{-1}\circ G_{\theta}'\circ G_{\phi}'\circ d=\Theta(G_{\theta}')\circ\Theta(G_{\phi}').
\]

Therefore $G(CF)$ is isomorphic to $G(OF)$.
\proofend

Since $G'$ and $G$ coincide for circles, we will from now on drop
the primes.

We fix one of the conics of the family $CF$ to be $(p+x)^{2}+y^{2}=(1+px)^{2}$. Such a conic has eccentricity $p$, though the case $p=1$ is degenerate, being a line rather than a parabola. A calculation shows that the conics of the family have equation

\begin{equation}  \label{eq:pencil}
(p+x)^{2}+y^{2}=(1+px)^{2}t
\end{equation}

for $t\in[0,\infty)$ as the parameter. We denote this pencil of conics
by $CF$.

\begin{proposition}
The pencil of conics $CF$ shares focus $F=(-p,0)$ and foliates $\mathbb{R}^{2}\setminus\{(x,y):x=-\frac{1}{p}\}$.
\end{proposition}

\proof
It is clear from the duality transformation that the conics share
focus $(-p,0)$. Two conics with differing values of $t$ are clearly
disjoint. Let $(x,y)\in\mathbb{R}^{2}\setminus\{(x,y):x=-\frac{1}{p}\}$.
If $x\neq-\frac{1}{p}$ then the ratio 
\[
\frac{(p+x)^{2}+y^{2}}{(1+px)^{2}}
\]

is defined. Since it is nonnegative, there is some $t$ which equals
to this ratio. If $x=-\frac{1}{p}$, then the equation 

\[
(p-\frac{1}{p})^{2}+y^{2}=0.
\]

has no solution since $p\neq\pm1$.
\proofend

\begin{figure}[hbtp]
\centering
\frame{\includegraphics[width=2.4in]{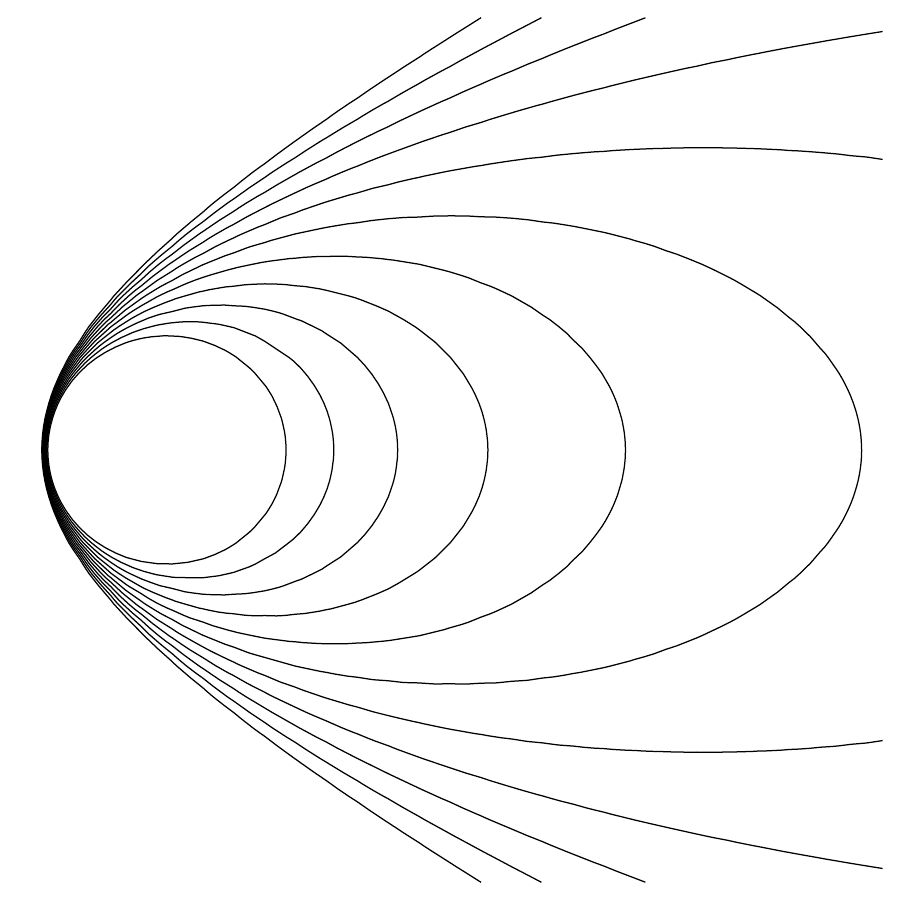}}
\frame{\includegraphics[width=2.4in]{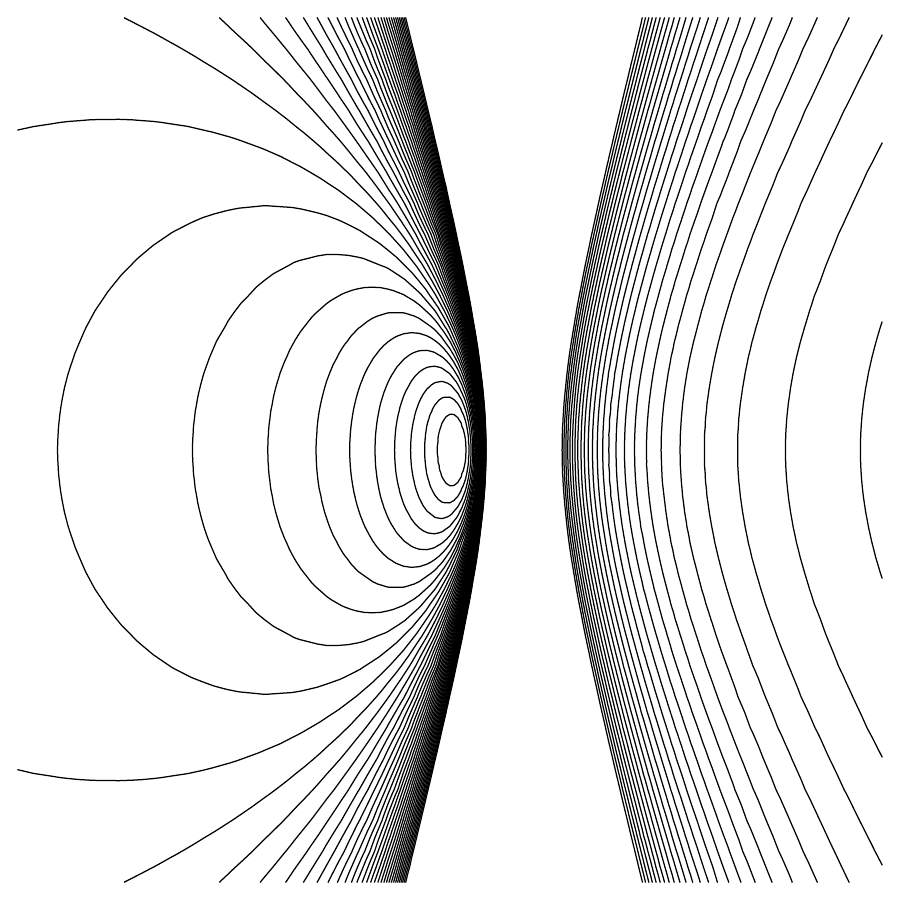}}
\caption{A portion of the pencil of conics $CF$ when $p=\frac{3}{4}$ and when $p=\sqrt{2}$. }
\end{figure}

Note that in the special case that $p=0$ we obtain a family of concentric
circles.

It is also important to note that the {\it discrete conics are precisely those polygons for
which there is an action of the group}. This action is obtained by taking $X=V_i$ and $Y=V_{i+k}$ to be vertices of the discrete conic. Then $\angle XFY=k \theta$.  
This is unique to discrete conics. Indeed, if $\tilde{D}$ is a polygon lying on a conic $C$ of the family $CF$ which is not a discrete conic, then at least one pair
of consecutive vertices do not form the same angle with respect to
the focus as the other vertices. Therefore the image of $\tilde{D}$,
call it $G(\tilde{D})$, will not lie on a single conic of $CF$.
Consequently, the operation $G^{2}(\tilde{D})$ will not be well-defined.

\section{The Discrete Negative Pedal Construction}

For a given curve $C$ in the plane and a given fixed point $P$, called the pedal point, the {\it pedal curve} of $C$ is the locus of points $X$ such that $PX$ is perpendicular to a tangent to the curve passing through $X$. 

The {\it negative pedal curve} is the inverse of the pedal curve. More precisely, the negative pedal curve of $C$ with respect to $P$ is the envelope of lines $XP$ for $X$ lying on the given curve. 
The negative pedal curve of a pedal curve with the same pedal point is the original curve.

For a conic $C$, if the pedal point is the focus $F$, then the pedal curve is a circle (a line in the case of the parabola). Conversely, given a circle $o$, if $P$ does not lie on $o$, then the negative pedal is a conic with $P$ as a focus. 

In \cite{Discrete Parabolas}, we used a discretized version of the negative pedal construction to construct discrete analogs of parabolas.  Our construction for the general conic is analogous except  {\it our construction before prioritized distance over angle}. For a general conic, no such choice is possible.
 We recall the previous construction. Let $P$ be a point,
$L$ a line and $X_{1},X_{2},\ldots,X_{n}$ points on $L$ with $|X_{i}X_{i+1}|=\Delta$
for $i=1,2,\ldots,n-1$. Let $L_{i}$ be the line orthogonal to $PX_{i}$
and passing through $X_{i}$, $i=1,2,\ldots,n-1$. The vertices of the
discrete parabola are $V_{i}=L_{i}\cap L_{i+1}$, $i=1,2,\ldots,n-1$. Note that distances are equal rather than angles. This gives an explanation as to why discrete parabolas turned out to be good metric approximators of the parabola, e.g., being the optimal piecewise-linear approximations to the parabola.

\begin{figure}[hbtp]
\centering
\includegraphics[width=2.4in]{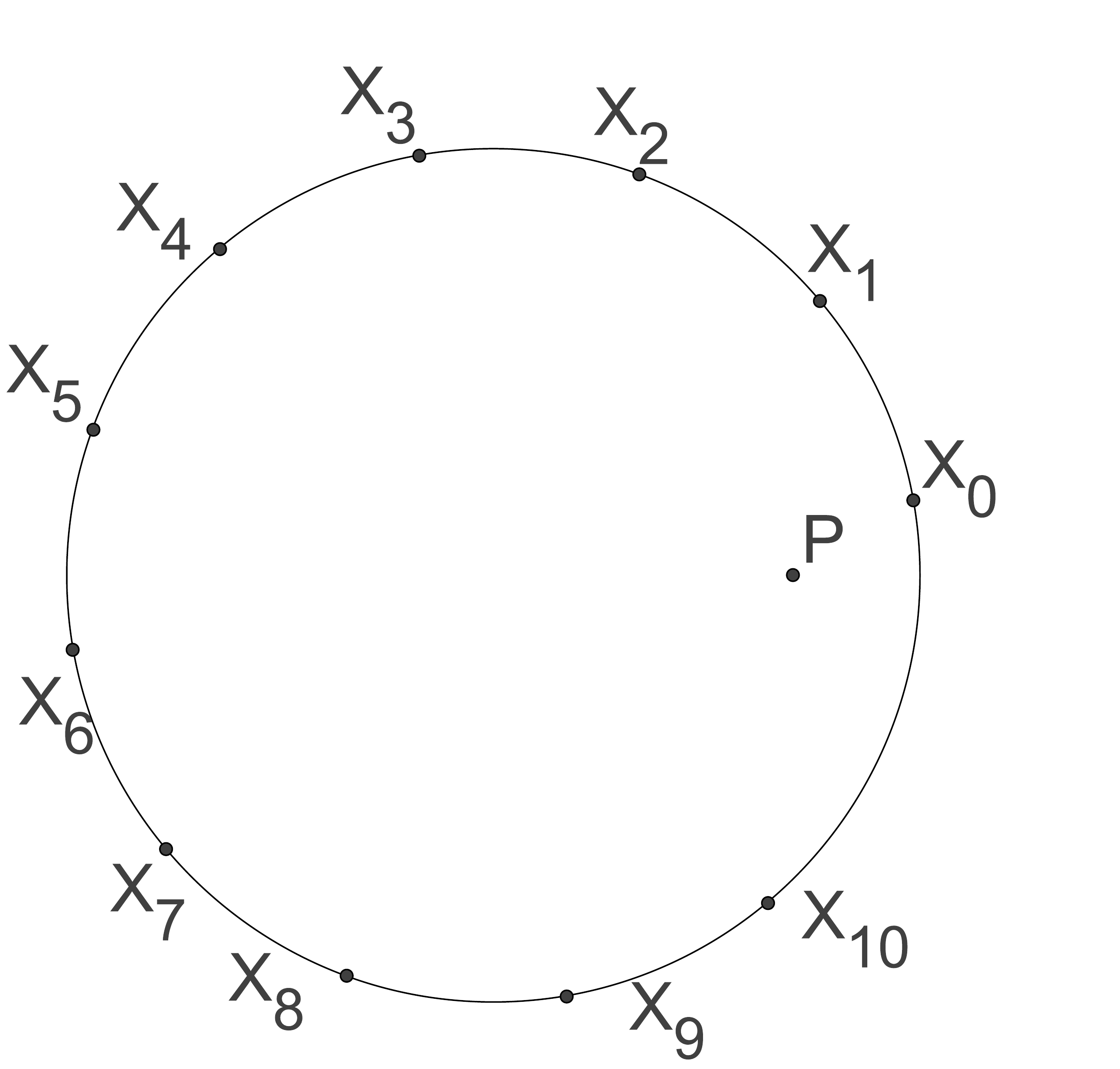}
\includegraphics[width=2.4in]{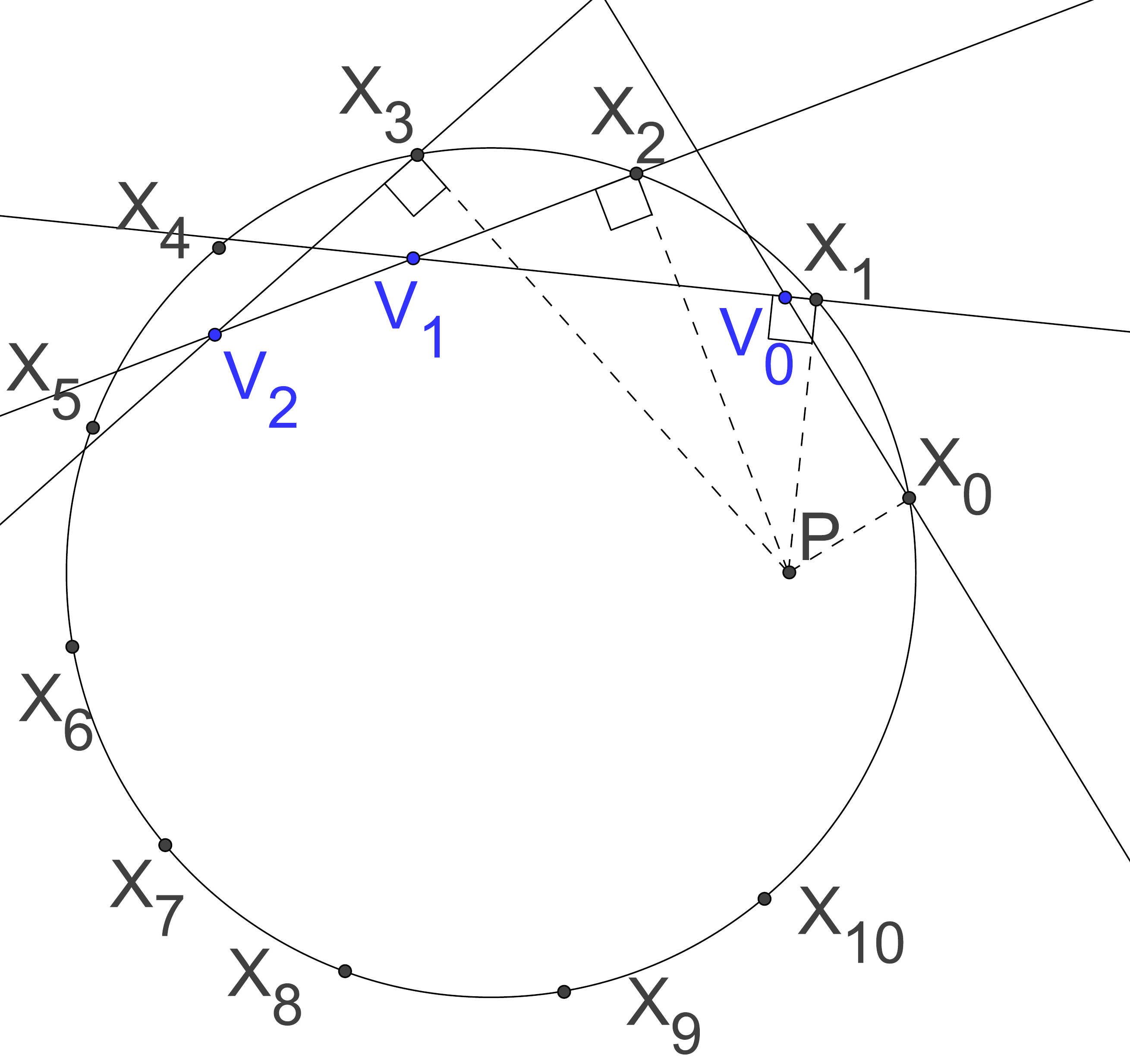}
\caption{\label{fig:Construction of Discrete Conic}Construction of a discrete
conic $V_{0},V_{1}V_{2}\vdots V_{n}.$ via the discrete negative pedal construction.}
\end{figure}

Now, we
replace the line $L$ by a circle. Let $P=(p,0)$ with $p\neq\pm1$
and let $X_{j}=(\cos(j\theta+\phi),\sin(j\theta+\phi))$, $j=1,2,\ldots,n$,
be equally spaced points on the unit circle. Let $L_{j}$ be the line
orthogonal to line $PX_{j}$ and passing through $X_{j}$. The vertices
of what we will later show to be a discrete conic are $V_{j+1}=L_{j}\cap L_{j+1}$, $j=0,1,2,\ldots,n-1$.
(see Figure \ref{fig:Construction of Discrete Conic}). 

A few remarks are in order. The cases $p=\pm1$ are degenerate for
our construction, so we do not consider them. We also assume that
the pedal point $P$ is on the $x$-axis because the general case can be obtained via rotation. 

We will show that the $V_{i}$'s lie on a conic independent of $\phi$.
We will also show that if $P'=(-p,0)$, then $\angle V_{i}P'V_{i+1}=\theta$
- a surprising fact considering that no mention of $P'$ is made in
the construction. 

In the limit  $\theta\rightarrow0$, we obtain the negative pedal construction with respect to a circle. Thus the negative pedal is the envelope of a conic $C_{L}$.  We call this conic the \textit{limiting conic.}
The sides of the discrete conic $V_{0}V_{1}\cdots V_{n-1}$ are
tangent to $C_{L}$, since they are a subset of the enveloping lines.

The following easy lemma will allow us to understand $C_{L}$ better.
\begin{lemma}

\label{lem:equal distances w.r.t circle center}Let $X$ and $Y$
be points on a circle $o$ centered at $O$. Let $B$ be line $XY$,
$L_{X}$ and $L_{Y}$ be the lines perpendicular to $B$ and passing
through $X$ and $Y$, respectively. Consider an arbitrary line $L$
passing through $O$. Set $P_{1}=L\cap L_{x}$ and $P_{2}=L\cap L_{Y}$.
Then the distance between $O$ and $P_{1}$ is the same as the distance
between $O$ and $P_{2}$.
\end{lemma}
\proof

The result clearly holds when $L\parallel B$. It is easy to see that
rotating $L$ preserves equality. 
\proofend

\begin{figure}[hbtp]
\centering
\includegraphics[width=2.4in]{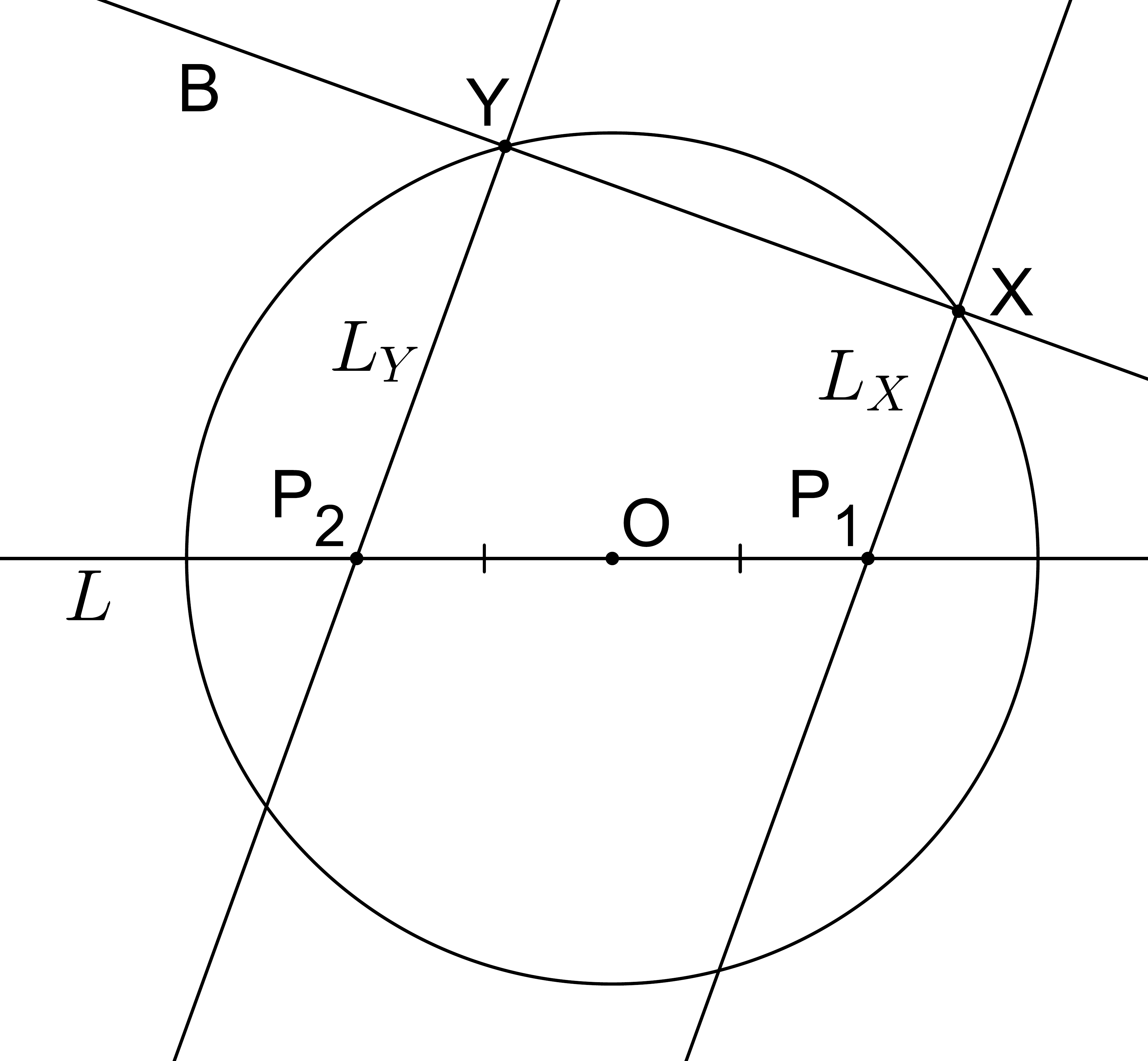}
\caption{Lemma \ref{lem:equal distances w.r.t circle center}.}
\end{figure}

\begin{theorem}

The focii of the limiting conic $C_{L}$ are $P=(p,0)$ and $P'=(-p,0)$.
The equation of $C_{L}$ is $x^{2}+\frac{y^{2}}{1-p^{2}}=1$.
\end{theorem}
\proof

We argue that the negative pedal of $P'$ with respect to the unit
circle is also $C_{L}$ so that $P'$ must be the other focus of $C_{L}$.
We do this by showing that the set of lines $L_{i}$ in the limit
$\theta\rightarrow0$ produced from $P$ is the same as those produced
from $P'$. This would imply that the negative pedal from $P'$ has
the same envelope of lines as that of $P$, so that its negative pedal
is $C_{L}$.

Let $X\in\mathbb{S}^{1}$. Let $B$ be the line orthogonal to $PX$
and passing through $X$. If $B$ is tangent to $\mathbb{S}^{1}$
then $X$ is on the $x$-axis, so it is also part of the envelope
corresponding to $P'$. Otherwise, $B$ intersects $\mathbb{S}^{1}$
in one more point - call it $Y$. Applying Lemma \ref{lem:equal distances w.r.t circle center},
we conclude that $B$ is also part of the envelope corresponding to
$P'$, so that $P$ and $P'$ produce the same envelope. It follows
that the focii of $C_{L}$ are $P$ and $P'$.

This implies that $C_{L}$ is of the form
\[
\frac{x^{2}}{\lambda}+\frac{y^{2}}{\lambda-p^{2}}=1.
\]

Since the line $x=1$ is part of the envelope, it must be tangent
to $C_{L}$. Therefore $C_{L}$ passes through $(1,0)$. It follows
that $\lambda=1$.
\proofend

The following two Lemmas will culminate in the result that the $V_{i}$'s
lie on a conic independent of $\phi$. The first may be viewed as
a simple result in mathematical billiards. The second follows from
the first, and is a possibly new property of the pedal circle with respect to its
conic.

\begin{figure}[hbtp]
\centering
\includegraphics[width=2.4in]{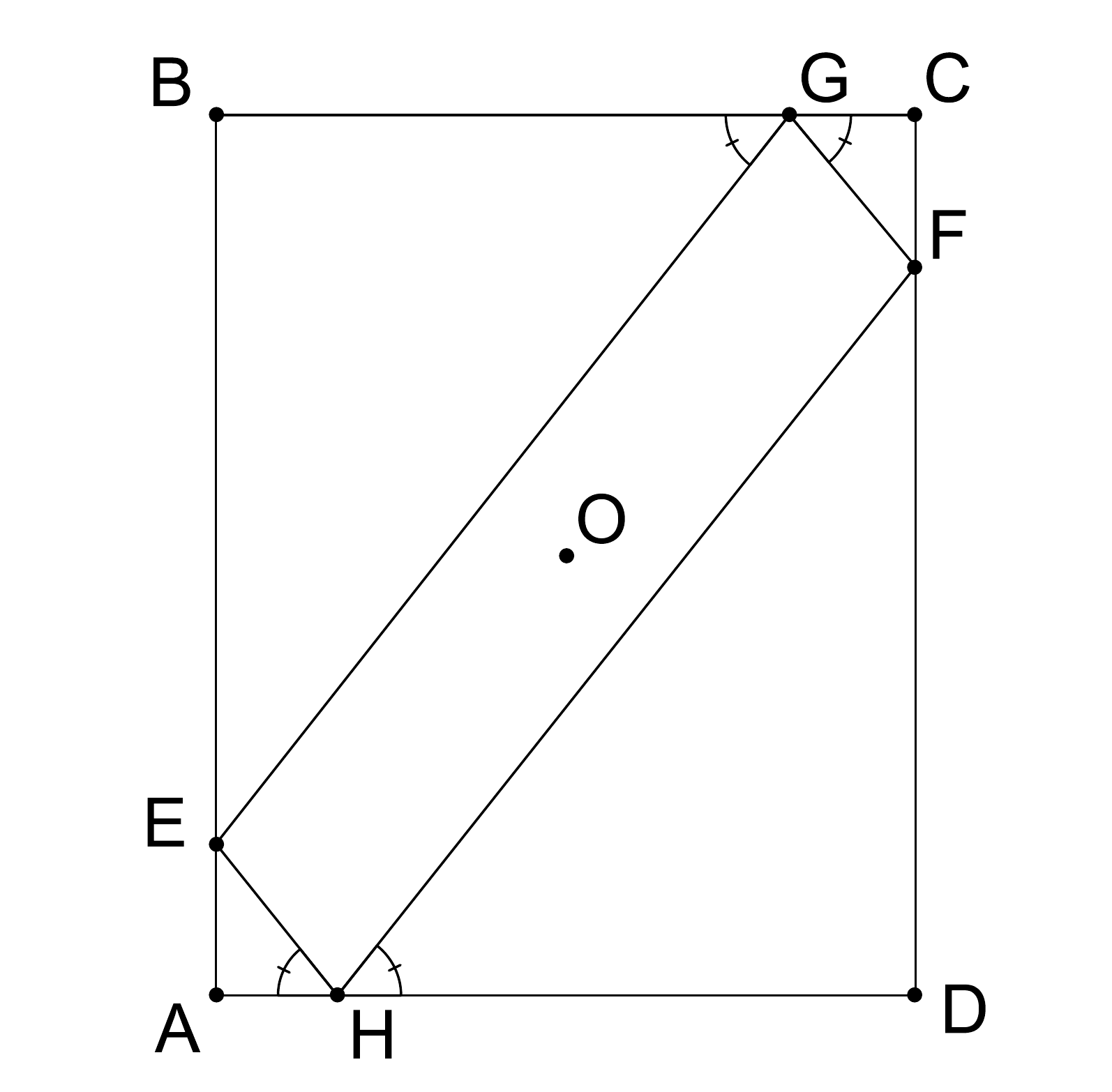}
\caption{Lemma \ref{lem:rectangle billiard}.}
\end{figure}

\begin{lemma}

\label{lem:rectangle billiard}Let $ABCD$ be a rectangle with center
$O$. Take points $E\in AB$ and $F\in CD$ such that $E,F,O$ are
collinear. Let $G\in BC$ and $H\in DA$ satisfy $\angle EGB=\angle FGC$
and $\angle EHA=\angle FHD$. Then

1. $E,G,F,H$ is a $4$-periodic billiard trajectory.

2. $\angle EGB=\angle FGC=\angle EHA=\angle FHD$ and $\angle AEH=\angle BEG=\angle CFG=\angle DFH$.

3. $G,H,O$ are collinear.

4. $EGFH$ is a parallelogram with sides parallel to the diagonals
of $ABCD$.
\end{lemma}
\proof

The equality $\angle EGB=\angle FGC=\angle EHA=\angle FHD$ and (3)
follow by odd symmetry about $O$. (1) and (2) then follow by a simple
angle-count. 

To prove (4), observe that $(1)$ implies that $EFGH$ is a parallelogram.
Let $\theta=\angle BEG$, $x=BE$ and $y=BG$. Triangle $\triangle EBG$
implies that
\[
\tan\theta=\frac{y}{x}.
\]
 Similarly, triangle $\triangle EAH$ implies that 
\[
\tan\theta=\frac{|BC|-y}{|AB|-x},
\]

so that $\triangle EBG$ is similar to $\triangle ABC$. 
\proofend

\begin{figure}[hbtp]
\centering
\includegraphics[width=2.4in]{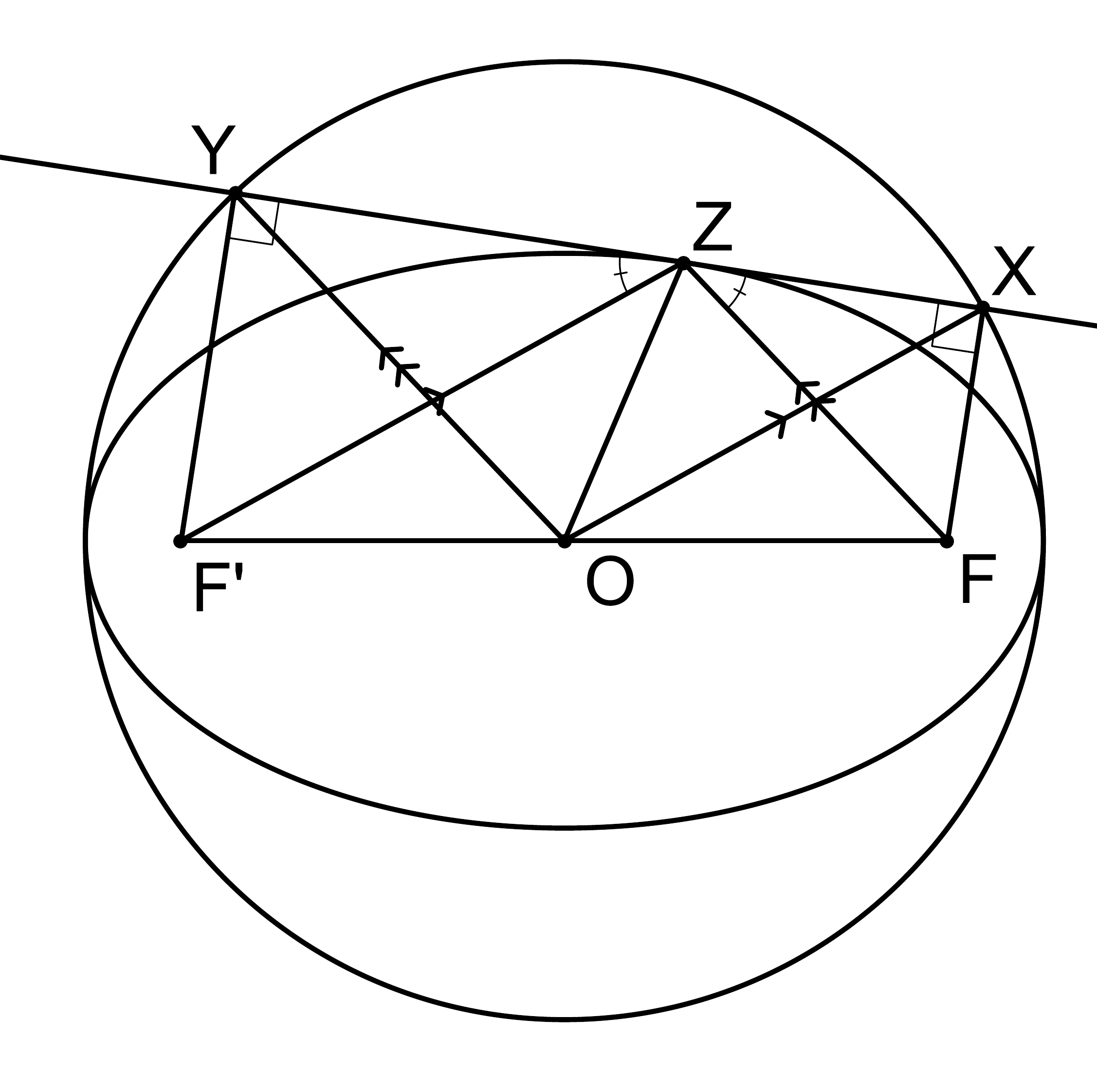}
\caption{Lemma \ref{lem:parallelism in conic and pedal circle}.}
\end{figure}

\begin{lemma}

\label{lem:parallelism in conic and pedal circle}Let $C$ be a conic
with focii $F$ and $F'$, let $o$ be the pedal circle with respect
to $F$ and call its center $O$. Consider any point $X$ on $o$.
Set $L$ to be the line through $X$ which is orthogonal to $FX$
and let its intersection with $o$ be $Y$ and its intersection with
$C$ be $Z$. Then $F'Z\parallel OX$ and $FZ\parallel OY$.
\end{lemma}
\proof

We extend $XF$ to intersect $o$ at $X'$ and $YF'$ to intersect
$o$ at $Y'$. Then $XYY'X'$ is a rectangle with center $o$. Call
the intersection of $X'Y'$ with $C$ $Z'$. We apply Lemma \ref{lem:rectangle billiard}
with $F$ and $F'$ being $E$ and $F$ from the Lemma and $Z$ and
$Z'$ being $G$ and $H$. 
\proofend

We are now ready to consider the construction at the beginning
of the section. Let $M_{j}=C_{L}\cap L_{j}$
for each $j$ and set $M=M_{1},M_{2} \cdots M_{n}$. A consequence of Lemma \ref{lem:parallelism in conic and pedal circle}
is the following result, which shows that $D=V_{1},V_{2} \dots $ and $M$ are discrete conics. 
\begin{theorem}

\label{thm:vertices lie on a conic}The intersections $V_{j}=L_{j}\cap L_{j+1}$
lie on a conic independent of $\phi$ and having focus $P'$ and form a discrete conic of angle $\theta$. Moreover, $D=G_{\theta}(M)$.
\end{theorem}
\proof
We argue that $\angle M_{j}P'M_{j+1}=\theta$. Indeed, line $P'M_{j}$
is parallel to $OX_{j}$ by Lemma \ref{lem:parallelism in conic and pedal circle}.
Since $\angle X_{j}OX_{j+1}=\theta$, the result follows. 

Recall that if tangents at $X$ and $Y$ to a conic intersect at $Z$
and $F$ is a focus of the conic then $ZF$ bisects angle $\angle XFY$
\cite{Akopyan Geometry of Conics}. It follows that $\angle V_{j}P'M_{j}=\frac{\theta}{2}$
and $\angle V_{j}P'V_{j+1}=\theta$. Using Theorem \ref{thm:conic locus}
with this last fact, we conclude that the $V_{j}$'s lie on a conic
independent of $\phi$ which shares focus $P'$ with $C_{L}$.
\proofend

\begin{figure}[hbtp]
\centering
\includegraphics[width=2.4in]{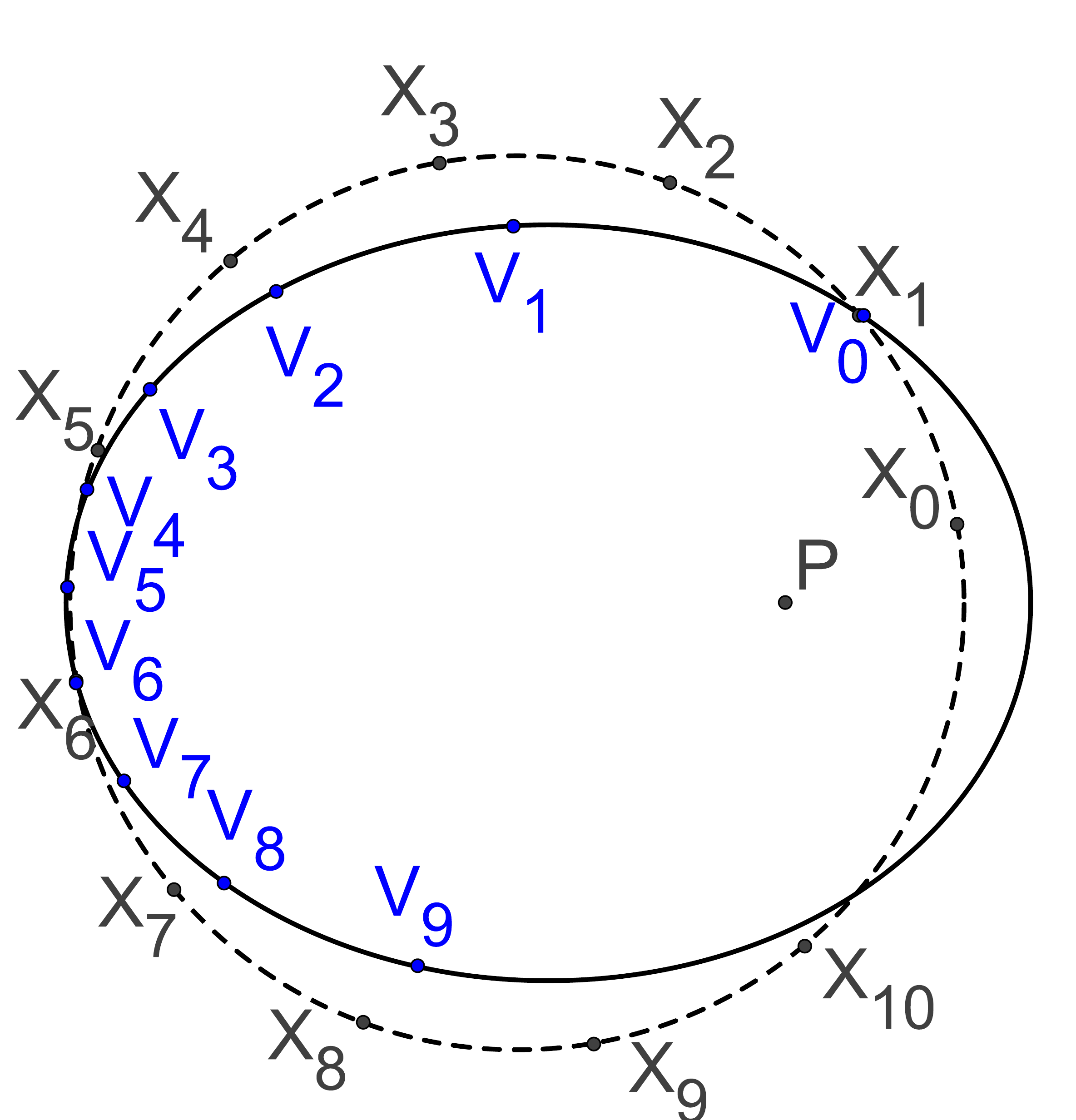}
\includegraphics[width=2.4in]{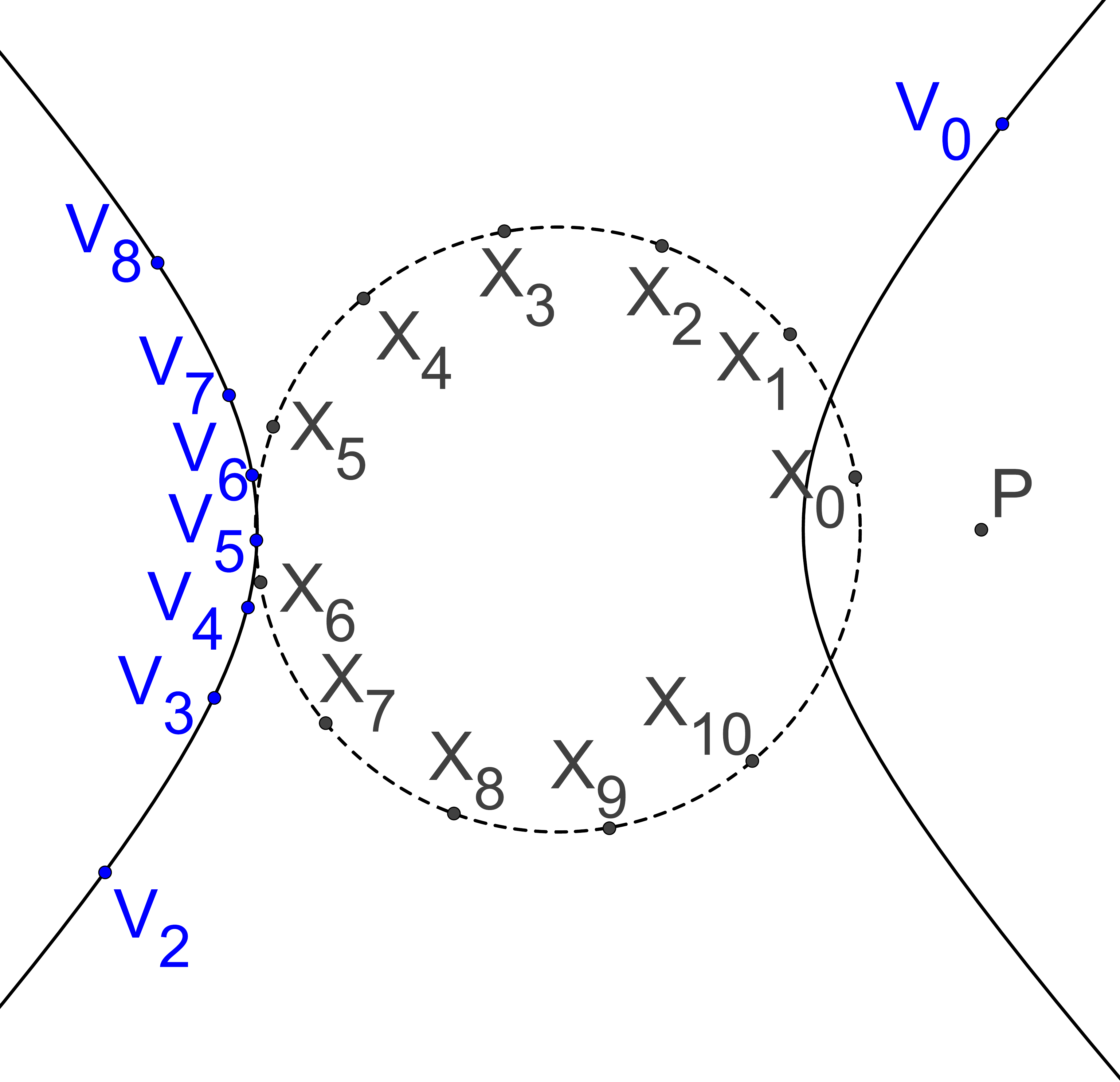}
\caption{Two different values of $p$. By Theorem \ref{thm:vertices lie on a conic},
the $V_{i}$'s lie on a conic independent of the phase of $X_{0},X_{1},...$ and form a discrete conic.}
\end{figure}

Let $C$ have equation $x^{2}+\frac{y^{2}}{1-p^{2}}=1$. A calculation
shows that for a discrete conic with parameter $\theta$ and phase
angle $\phi$ has vertices with coordinates 
\[
V_{j}=(\frac{p-\cos((j-1)\theta+\phi)}{p\cos((j-1)\theta+\phi)-1},\frac{(p^{2}-1)\sin((j-1)\theta+\phi)}{p\cos((j-1)\theta+\phi)-1}).
\]

\section{Proofs of Theorems from Section \ref{intro}} \label{Proofs}

\noindent{\bf Proof of Theorem  \ref{Thm: tangency points projectively regular}.} Let $D=V_1 V_2 \cdots V_{n}$ and let $T_i$ denote the tangent at $V_i$ to the conic $C$. Applying the duality transformation $d$ as in the proof of Theorem \ref{thm:conic locus}, we see that each $d(V_i)$ is a tangent line to circle $d(C)$, and $d(T_i)$ becomes the the tangency point of $d(V_i)$ on $d(C)$. Let $p_i$ denote the projection of $F$ into $d(V_i)$. Since $\angle V_i F V_{i+1}=\theta$ for each $i$, we have $\angle p_i F p_{i+1}=\theta$. This implies that the angle formed by consecutive tangents is $\pi-\theta$, and in particular, is always equal. This shows that the vertices of $D$ map to the sides of a regular polygon. The sides of $D$ map to the vertices of this regular polygon.
\proofend

\noindent{\bf Proof of Theorem \ref{Thm: Poncelet Polygon}.} This follows from Theorem \ref{Thm: tangency points projectively regular} but we can provide an additional proof. Given a discrete conic $D$ on a conic $C$ with angle parameter $\theta$, the discrete conic $D'=H_{\theta}(D)$ lies on a conic $C'$ sharing a focus with $C$ and satisfying $D=G_{\theta}(D')$. This means that the vertices of $D$ are intersections of tangents to $C'$ at  consecutive pairs of  vertices of $D'$. In particular, the sides of $D$ are tangent to $C'$.
\proofend

\noindent{\bf Proof of Theorem \ref{thm:Diagonals coincide at P'}.} This too follows from Theorem \ref{Thm: tangency points projectively regular} but we provide another proof. Assume first that $n$ is even. Using the fact that $\angle V_{j}FV_{j+1}=\theta$
we see that $\angle V_{j}FV_{j+\frac{n}{2}}=\theta\frac{n}{2}$. Since
$n=\frac{2\pi}{\theta}$, we conclude that $V_{j},F$ and $V_{j+\frac{n}{2}}$
are collinear. 

In the odd case, we have $\angle V_{i}FM_{i}=\frac{\theta}{2}$ and
$\angle M_{i}FM_{i+1}=\theta$. Therefore $\angle V_{i}FM_{i+\frac{n-1}{2}}=\frac{\theta}{2}+\frac{n-1}{2}\theta=\frac{n}{2}\theta$,
so that $V_{j},F$ and $M_{j+\frac{n-1}{2}}$ are collinear.
\proofend

\noindent{\bf Proof of Theorem \ref{thm:Reflective Property}.} Assume first that the number of sides $n$ of $D$ is even. Since
$r$ traverses $F'\rightarrow S_{j}\rightarrow F$, and $S_{j}$ is
a line tangent to $C_{L}$, the reflection of $r$ in $S_{j}$ passes
through $M_{j}$ and $F$. Theorem \ref{thm:Diagonals coincide at P'}
implies that $M_{j}$, $F$ and $M_{j+\frac{n}{2}}$ are collinear,
since $M_{1}M_{2}\cdot\cdot\cdot M_{n}$ is also a discrete conic
with parameter $\theta$ and focus $F$. Therefore $r$ must subsequently
reflect in side $S_{j+\frac{n}{2}}$ at $M_{j+\frac{n}{2}}$ and pass
through $F'$.

Now assume that $n$ is odd. Then Theorem \ref{thm:Diagonals coincide at P'}
implies that $r$ passes through $V_{j-\frac{n-1}{2}}=V_{j+\frac{n-1}{2}}$.
\proofend

\noindent{\bf Proof of Theorem \ref{thm:isogonal property}.}
The lines $S_{i}$ and $S_{j}$ are tangent to $C_{L}$ at $M_{i}$
and $M_{j}$. Since $F$ is a focus of $C_{L}$ and $Z$ is the intersection
of the tangents, $FZ$ bisects angle $\angle M_{i}FM_{j}$. Since
$\angle M_{i}FV_{i}=\angle M_{j}FV_{j+1}$, we see that $FZ$ bisects
angle $\angle V_{i}FV_{j+1}$. Similar reasoning shows that $FZ$
bisects angle $\angle V_{i+1}FV_{j}$. 

Finally, since $S_{i}$ and $S_{j}$ are tangent to $C_{L}$, the
isogonal property of a conic implies that $FZ$ and $F'Z$ are isogonal
with respect to $S_{i}$ and $S_{j}$.
\proofend

\noindent{\bf Proof of Theorem \ref{thm:Poncelet Grid}.} Analogous to the proof of Theorem \ref{Thm: Poncelet Polygon}.
\proofend

\noindent{\bf Proof of Theorem \ref{thm: pascal}.} The intersections of opposite lines are the vertices of the discrete conic $G_{n\theta}(H_{\theta}(D))$, implying that they are a discrete conic.
In particular, the intersection of sides $S_i S_{i+n}$ maps to the line $d(S_i) d(S_{i+n})$ under duality. Since $d(S_1) d(S_2) \cdots d(S_{2n})$ form a regular polygon, the dual of the intersection of opposite sides are the $2n$ diagonals of the regular polygon (each diagonal is counted twice). These pass through the center $O$ of the cirumscribing circle of this regular polygon. Note that $F'$ is the image of the center of $O$. Applying duality once more, we see that these diagonals map to points lying on the dual of $O$, which is a line orthogonal to $FF'$. 
\proofend

\bigskip
{\bf Acknowledgments}. The author is grateful to S. Tabachnikov, R. Schwartz and J. Sondow for their helpful suggestions.

\end{document}